%% file: SDE.tex
\documentclass[onecolumn,a4paper,10pt]{siamart190516}

\input{ex_shared}

\ifpdf
\hypersetup{
  pdftitle={Polynomial propagation of moments in stochastic differential equations},
  pdfauthor={A. López Yela and J. Míguez}
}
\fi

\begin{document}

\maketitle

\begin{abstract}
\cblack{
We address the problem of approximating the moments of the solution, $\tb{X}(t)$, of an It\^o stochastic differential equation (SDE) with drift and a diffusion terms over a time-grid $t_0, t_1, \ldots, t_n$. In particular, we assume an explicit numerical scheme for the generation of sample paths $\hat{\tb{X}}(t_0), \hat{\tb{X}}(t_1), \ldots, \hat{\tb{X}}(t_n), \ldots$ and then obtain recursive equations that yield any desired non-central moment of $\hat{\tb{X}}(t_n)$ as a function of the initial condition $\hat{\tb{X}}(t_0) = \tb{X}_0$. The core of the methodology is the decomposition of the numerical solution $\tb{X}(t_n)$ into a ``central part'' and an``effective noise'' term. The central term is computed deterministically from the ordinary differential equation (ODE) that results from eliminating the diffusion term in the SDE, while the effective noise accounts for the stochastic deviation from the numerical solution of the ODE. For simplicity, we describe the proposed methodology based on an Euler-Maruyama integrator, but other explicit numerical schemes can be exploited in the same way. We also apply the moment approximations to construct estimates of the 1-dimensional marginal probability density functions of $\hat{\tb{X}}(t_n)$ based on a Gram-Charlier expansion. Both for the approximation of moments and 1-dimensional densities, we describe how to handle the cases in which the initial condition is fixed (i.e., $\tb{X}_0 = \tb{x}_0$ for some known $\tb{x_0}$) or random. In the latter case, we resort to polynomial chaos expansion (PCE) schemes to approximate the target moments. The methodology has been inspired by the PCE and differential algebra (DA) methods used for uncertainty propagation in astrodynamics problems. Hence, we illustrate its application for the quantification of uncertainty in a 2-dimensional Keplerian orbit perturbed by a Wiener noise process. 
}
\end{abstract}



\section{Introduction}\label{introduction}

\cblack{
Let us consider the stochastic differential equation (SDE) in It\^o form \cite{Oksendal13}
\begin{equation}\label{dyn_sde}
	\left\{\begin{matrix}
		\diff \tb{X}(t)&=&\tb{u}(\tb{X},t)\diff t+\tb{G}(\tb{X},t)\diff\tb{W}(t)\\
	 	\tb{X}(0) &= &\hspace{-3.45cm}\tb{X}_0,
	\end{matrix}\right.
\end{equation}
where $\tb{X}(t)$, $t \ge 0$, is a real $v$-dimensional random process representing the solution of the SDE, $\tb{X}_0$ is a real $v$-dimensional random variable that describes the initial condition of the process, functions $\tb{u}:\mathds{R}^v\times\mathds{R}_0^{+}\rightarrow\mathds{R}^v$ and $\tb{G}:\mathds{R}^v\times\mathds{R}_0^{+}\rightarrow\mathds{R}^{v\times d}$ are the the drift coefficient and the diffusion coefficient, respectively, and $\tb{W}(t)$ is a  $d$-dimensional stochastic process with independent increments.
}

\cblack{
When $\tb{W}(t)$ is assumed to be a Wiener process and the drift and diffusion coefficients satisfy some standard differentiability assumptions, it can be shown that the solution $\tb{X}(t)$ to Eq. \eqref{dyn_sde} can be characterised by a time-varying probability density function (pdf) that we denote as $f_{\tb{X}}$ and satisfies the Fokker--Planck equation \cite{Ri89}
\begin{equation}\label{FPeq}
\frac{\partial f_{\tb{X}}(\tb{x},t)}{\partial t}+\sum_{k=1}^v \frac{\partial}{\partial \tb{x}^{(k)}}[u^{(k)}\!(\tb{x},t) f_{\tb{X}}(\tb{x},t)]-\frac{1}{2}\sum_{k=1}^v\sum_{j=1}^v\frac{\partial^2}{\partial \tb{x}^{(k)}\partial \tb{x}^{(j)}}[D^{(kj)}(\tb{x},t) f_{\tb{X}}(\tb{x},t)]=0
\end{equation}
with initial condition $f_{\tb{X}}(\tb{x},0) = f_{\tb{X}_0}(\tb{x})$, where $u^{(k)}$, $k=1, \ldots, v$, are the components of the drift coefficient $\tb{u}(\tb{X},t)$ in Eq. \eqref{dyn_sde} and $D^{(kj)}(\tb{x},t)$ is the entry in the $k$-th row and $j$-th column of the diffusion tensor $\tb{D}(\tb{x},t)=\tb{G}(\tb{x},t)\tb{G}(\tb{x},t)^\top$. In principle, we may completely characterise the solution of Eq. \eqref{dyn_sde} by solving the partial differential equation (PDE) \eqref{FPeq}. However, this cannot be done exactly except for special (simple) cases \cite{Ri89}. On the other hand, the computational cost of numerical schemes for PDEs, based on finite differences \cite{Sm85} or finite elements \cite{Br08,LY17}, quickly becomes prohibitive as the dimensions $v$ and $d$ increase.
}   

\cblack{
Because of the difficulties in solving the Fokker-Planck equation \eqref{FPeq}, most authors have focused on the study of time-discretisation numerical schemes to simulate realisations of the random process $\tb{X}(t)$. Such schemes are extensions of classical algorithms for the numerical solution of ordinary differential equations (ODEs) and they include the classical Euler-Maruyama, Milstein or stochastic Runge-Kutta methods \cite{Gard88,Kloeden92}, as well as their implicit and semi-implicit variants \cite{Tian01,Mao13,Yao18}. When the noise process $\tb{W}(t)$ is Wiener, the convergence and stability of these numerical algorithms can be studied using a variety of techniques \cite{Gard88,Kloeden07,Hutzenthaler12}, although Taylor approximations have become the standard approach in the past years \cite{Kloeden92}. Let us remark, however, the fundamental difference between simulating a realisation $\tb{X}(t)=\tb{x}(t)$ for a discrete-time grid, $t \in \{ t_0, t_1, \ldots, t_N \}$, and the probabilistic characterisation that would be obtained by computing the pdf's $f_{\tb{X}}(\tb{x},t_i)$, even if just approximately. While one can certainly generate many trajectories $\tb{X}(t)=\tb{x}_i(t)$, $i=1, ..., N$, in order to construct Monte Carlo estimators over the grid $t \in \{ t_0, t_1, \ldots, t_N \}$, the computational cost of such an approach becomes intractable, again, as the dimension $v$ of the process increases.
}

\cblack{
In this paper we introduce a new approach to the probabilistic characterisation of the solution $\tb{X}(t)$ to the SDE \eqref{dyn_sde}. Choose a time grid $\{t_0, \ldots, t_n, \ldots\}$, an initial condition $\tb{X}_0=\tb{x}_0$ and let $\hat{\tb{X}}_n$ be the random sequence generated by the Euler-Maruyama scheme applied to the SDE \eqref{dyn_sde}. The proposed method builds upon: 
\begin{itemize}
\item[(a)] The classical Euler scheme applied to the ordinary differential equation (ODE) $\dot{\tb{x}} = \tb{u}(\tb{x},t)$ with initial condition $\tb{x}_0$, that yields a deterministic sequence $\hat{\tb{x}}^C_n \simeq \tb{x}(t_n)$, $n=0,1,...$. We refer to this sequence as the \textit{central} part of $\hat{\tb{X}}_n$. 
\item[(b)] The construction of an \textit{effective noise} sequence, denoted $\Delta\hat{\tb{W}}_n$, that relates the central component and the Euler-Maruyama realisation as $\hat{\tb{X}}_n=\hat{\tb{x}}^C_n + \Delta \hat{\tb{W}}_n$.
\end{itemize} 
We show how the moments of the effective noise process $\Delta\hat{\tb{W}}_n$ can be approximated recursively using a polynomial (Taylor) expansion. The moments of the random vectors $\hat{\tb{X}}_n$ are then obtained in a straightforward way via the binomial theorem. When the initial condition $\tb{X}_0$ is random, the method can be combined in a straightforward way with a polynomial chaos expansion (PCE) scheme to account for the initial uncertainty. Finally, we also show how to approximate the marginal pdf of each component $\hat X_n^{(k)}$ in the vector $\hat{\tb{X}}_n = \left( \hat{X}_n^{(1)}, \ldots, \hat{X}_n^{(v)} \right)$ by combining its moment estimates with a Gram-Charlier expansion of type A. The practical performance of the proposed scheme is illustrated with two examples related to astrodynamics, namely the propagation of uncertainty for a Keplerian orbit in two dimensions perturbed by a Wiener process. 
}

\cblack{ 
While in this manuscript we have restricted the analysis to the Euler and Euler-Maruyama schemes for the sake of clarity, our arguments can be extended to other numerical algorithms.
}

\cblack{
The rest of the paper is organised as follows. In Section \ref{sec_algorithm}, we introduce the methodology and outline the recursive algorithms for the approximation of the moments of $\hat{\tb{X}}_n$ with both fixed ($\tb{x}_0$) and random ($\tb{X}_0$) initial condition, as well as the scheme to estimate the marginal pdf's of $\hat{X}_n^{(k)}$, $k=1, \ldots, v$, from the approximate moments. In Sections \ref{var_sol_noise} and \ref{pdf_sol_noise} we present the analysis that supports the proposed algorithms. In particular, in Section \ref{var_sol_noise} we establish the convergence of the estimates of the moments of the effective noise terms (when the order of their polynomial approximations increases), while in Section \ref{pdf_sol_noise} we provide conditions for the convergence of the Gram-Charlier expansion of the marginal pdf's. In Section \ref{Num_examples} we apply the proposed numerical schemes to the characterisation of the uncertainty in a 2-dimensional Keplerian orbit perturbed by a Wiener process. of an object in two dimensions. Finally, a discussion of the theoretical and numerical results is presented in Section \ref{sConclusions}. 
}


\section{The algorithm}\label{sec_algorithm}

\cblack{In this Section we introduce the proposed algorithms for the approximation of the moments and the 1-dimensional marginal pdf's of the solution of Eq. \eqref{dyn_sde} over a time grid. These schemes are the main contribution of the paper. We provide the general argument for their derivation and a summary aimed at facilitating their implementation, but postpone the proof of the key theoretical results to Sections \ref{var_sol_noise} and \ref{pdf_sol_noise} for clarity. We start with a brief summary of the key notation used in this section (and the rest of the paper).}

\subsection{Notation}

\cblack{Consider a probability space $(\Omega,\mF,\mbP)$, where $\Omega$ is the sample space, $\mF$ denotes a $\sigma$-algebra of subsets of $\Omega$ and $\mbP$ is a reference probability measure. We denote random variables (r.v.'s) and random processes (r.p.'s) on $(\Omega,\mF,\mbP)$ with capital letters, e.g., $X$ and $X(t)$, respectively, and use lower-case letters to indicate specific realisations. For example, $x$ is a realisation of the r.v. $X$ and $x(t)$ denotes a sample path of $X(t)$.}

\cblack{Vectors are denoted with bold-face letters while we use regular-face for scalars, e.g., $\tb{x}$ and $x$, respectively. For a vector $\tb{x}$, a superscript $^{(k)}$ indicates the $k$-th component of the vector, i.e., if $\tb{x}$ is a $v$-dimensional vector then $\tb{x}=\left( x^{(1)}, \ldots, x^{(v)} \right)$. A multi-index $\tb{r}=\left( r^{(1)}, \ldots, r^{(v)} \right)$ is a vector of non-negative integers, i.e., $x^{(i)} \in \mbN \cup \{0\}$ for every $i$. We define the following shorthands for typical operations on multi-indices:}
\cblack{
\beqa
|\tb{r}| &\coloneq& \sum_{k=1}^v\tb{r}^{(k)}, \nn \\
\tb{r}! &\coloneq& \prod_{k=1}^v r^{(k)}!, \nn \\
\frac{
	\partial^{|\tb{r}|}
}{
	\partial \tb{a}^{\tb{r}}
} &\coloneq& \frac{
	\partial^{|\tb{r}|}
}{
	\partial a^{(1)^{r^{(1)}}} \cdots~~\partial a^{(v)^{r^{(v)}}}
} \nn \\
\tb{a}^{\tb{r}} &\coloneq& \prod_{k=1}^v a^{r^{(k)}} \nn \\
\sum_{\tb{r}'=\tb{0}}^{\tb{r}} &\coloneq& \sum_{r'^{(1)}=0}^{r^{(1)}}\cdots \sum_{r'^{(v)}=0}^{r^{(v)}},\nn\\
\begin{pmatrix}
\tb{r}\\
\tb{r}'\\
\end{pmatrix} &\coloneq& \prod_{k=1}^v 
\begin{pmatrix}
r^{(k)}\\
r'^{(k)}\\
\end{pmatrix}. \nn
\eeqa
We adopt the convention $0!=1$ (hence, $(0,\ldots,0)!=1$ as well).
}

\subsection{Euler--Maruyama discretisation and the effective noise process}

\cblack{The discretisation of the SDE \eqref{dyn_sde} using the explicit Euler--Maruyama scheme yields:}
\begin{equation}\label{sto_integrator}
\hat{\tb{X}}_n=\hat{\tb{X}}_{n-1}+h\tb{u}\big(\hat{\tb{X}}_{n-1},t_{n-1}\big)+\tb{G}\big(\hat{\tb{X}}_{n-1},t_{n-1}\big)\Delta\tb{W}_{n},\qquad n=1,2,\ldots,
\end{equation}
\cblack{where $\hat{\tb{X}}_n\simeq\tb{X}(t_n)$ is the approximation of the solution at time $t_n$, with $t_n=t_0+nh$, the subscript $n$ denotes discrete time, $h$ is the step-size and $\Delta\tb{W}_{n}=\tb{W}(t_n)-\tb{W}(t_{n-1})$ is the increment of the r.p. $\tb{W}(t)$ in the interval $(t_{n-1},t_n)$. The key of the proposed method is to decompose the random sequence $\tb{X}_n$ into two parts: a \textit{central part}, that results from the integration of an ODE, and an \textit{effective noise} sequence that accounts for the randomness in $\tb{X}_n$. These two notions are explicitly introduced below.} 
\begin{definition}\label{sto_split}
\cblack{The random sequence in Eq. \eqref{sto_integrator} can be written as 
$$
\hat{\tb{X}}_n=\hat{\tb{x}}_n^C+\Delta\hat{\tb{W}}_n,
$$
where the deterministic sequence $\hat{\tb{x}}_n^C\simeq\tb{x}(t_n)$ is the central part that results from the explicit Euler integration of the ODE $\dot{\tb{x}} = \boldsymbol{u}(\tb{x},t)$ with a prescribed initial condition $\tb{x}^{C}(t_0)=\tb{x}_0$; specifically
\beq
\hat{\tb{x}}_n^C=\hat{\tb{x}}_{n-1}^C+h\tb{u}\big(\boldsymbol{\hat{x}}_{n-1}^C,t_{n-1}\big),\qquad n\in\mathbb{N},
\label{eqCentralPart}
\eeq
and $\Delta\hat{\tb{W}}_n = \hat{\tb{X}}_n - \hat{\tb{x}}_n^C$ is the effective noise r.p.
}
\end{definition}

\cblack{The central part is easily computed as in Eq. \eqref{eqCentralPart}. However, the characterisation of the effective noise is not straightforward. The gist of our approach is to perform a Taylor expansion of $\Delta\hat{\tb{W}}_n$ around $\hat{\tb{x}}_{n-1}^C$ at each time step. Such expansion is convenient because it naturally provides a probabilistic description of the effective noise (and, as a consequence, of the numerical solution $\hat{\tb{X}}_n$) and it can be carried out recursively over time.}


\subsection{Polynomial expansion of the effective noise}\label{pol_eff_noise_A}

\cblack{For the analysis of the effective noise process it is convenient to handle separately the uncertainty in Eq. \eqref{sto_integrator} due to the random initial condition $\tb{X}_0$ and the uncertainty due to the sequence of independent noise increments $\Delta\tb{W}_n$ (this separation is already implicit in the definition of the effective noise). Consequently, let us first assume that the initial condition is deterministic and fixed\footnote{We analyse the case with a random initial condition in Section \ref{PCE_subsec}.}, i.e., $\tb{X}_0=\tb{x}_0$. The probability distributions and statistical moments of the r.p.'s $\hat{\tb{X}}_n$ and $\Delta\hat{\tb{W}}_n$ can then be computed conditional on $\tb{X}_0=\tb{x}_0$. In particular, in \textbf{Section \ref{var_sol_noise}}, we prove that the polynomial expansion of order $N$ for the effective noise at time $n$ with initial condition $\tb{X}_0=\tb{x}_0$ can be recursively written as
\beqa
\Delta\hat{W}_{n,N}^{(k)}(\tb{x}_0) &=& \Delta\hat{W}_{n-1,N}^{(k)}(\tb{x}_0) \nn \\
&& +h\!\!\sum_{|\boldsymbol{r}|=1}^N\frac{1}{\boldsymbol{r}!}\frac{\partial^{|\boldsymbol{r}|}}{\partial \tb{x}_{n-1}^{\boldsymbol{r}}}u^{(k)}\big(\hat{\tb{x}}_{n-1}^C(\tb{x}_0),t_{n-1}\big)\Delta\hat{\tb{W}}_{n-1,N}^{\boldsymbol{r}}(\tb{x}_0) \nn \\
&&+\sum_{|\boldsymbol{r}|=0}^{N-1}\sum_{j=1}^d\frac{1}{\boldsymbol{r}!}\frac{\partial^{|\boldsymbol{r}|}}{\partial \tb{x}_{n-1}^{\boldsymbol{r}}}G^{(k,j)}\big(\hat{\tb{x}}_{n-1}^C(\tb{x}_0),t_{n-1}\big)\Delta W^{(j)}_{n}\Delta\hat{\tb{W}}_{n-1,N}^{\boldsymbol{r}}(\tb{x}_0),
\label{eff_noise_expansion}
\eeqa
for $k=1,\ldots,v$ and, hence, we denote $\Delta\hat{\tb{W}}_{n,N}(\tb{x}_0)=\big(\Delta\hat{W}_{n,N}^{(1)}(\tb{x}_0),\ldots,\Delta\hat{W}_{n,N}^{(v)}(\tb{x}_0)\big)$. Since we have assumed that the initial condition is fixed, then, $\Delta\hat{\tb{W}}_{0,N}(\tb{x}_0)=\textbf{0}$. Notice that the multi-index $\tb{r}$ in the summations is $v$-dimensional. The subscript $N$ in $\Delta\hat{\tb{W}}_{n,N}(\tb{x}_0)$ indicates that we construct a polynomial approximation of order $N$ with no remainder term.}

\cblack{From Eq.\,\eqref{eff_noise_expansion}, it is straightforward to obtain the expansion of $\Delta\hat{\tb{W}}_{n,N}^{\tb{r}}(\tb{x}_0)$ (using combinatorics) for any $v$-dimensional multi-index $\tb{r}$ such that  $|\tb{r}|>1$. In particular, the conditional moments of the effective noise truncated to order $N$ can be written as
\begin{equation}\label{cond_moms_eff_noise}
\mathds{E}\big[\Delta\hat{\tb{W}}^{\tb{r}}_{n,N}(\tb{x}_0)\big]=\!\!\sum_{|\boldsymbol{s}|+|\boldsymbol{r'}|=1}^N \!\!\!a_{\tb{r},n,N}^{\tb{s},\boldsymbol{r'}}\big(\hat{\tb{x}}_{n-1}^C(\tb{x}_0),t_{n-1}\big)\mathds{E}\big[\Delta\tb{W}_{n}^{\,\boldsymbol{s}_{\vphantom{a_1}}}\big]\mathds{E}\big[\Delta\hat{\tb{W}}_{n-1,N}^{\boldsymbol{r'}}(\tb{x}_0)\big],
\end{equation}
where $a_{\tb{r},n,N}^{\tb{s},\boldsymbol{r'}}\big(\hat{\tb{x}}_{n-1}^C(\tb{x}_0),t_{n-1}\big)$ are the coefficients of the expansion obtained from the polynomial coefficients of Eq.~~\eqref{eff_noise_expansion}. Note that the multi-index $\boldsymbol{r'}$ is $v$-dimensional and the multi-index $\boldsymbol{s}$ is $d$-dimensional.}

\cblack{To obtain the identity \eqref{cond_moms_eff_noise}, we have assumed that the noise increments $\Delta\tb{W}_{n}$ form an independent random sequence, which implies that the effective noise $\Delta\hat{\tb{W}}_{m}$ is itself independent of $\Delta\tb{W}_{n}$ for every $m<n$.}

\cblack{Finally, using the binomial theorem, we arrive at a the formula of the conditional moments of $\hat{\tb{X}}_n$ given the initial condition $\tb{X}_0=\tb{x}_0$ in terms of the conditional moments of the effective noise in Eq.\,\eqref{cond_moms_eff_noise} and the central part,
\begin{equation}\label{cond_moms_sol}
\mathds{E}\big[\hat{\tb{X}}^{\boldsymbol{r}}_{n,N}(\tb{x}_0)\big]=\sum_{\boldsymbol{r'}=0}^{\boldsymbol{r}}\begin{pmatrix}
\boldsymbol{r}\\
 \boldsymbol{r'}
\end{pmatrix}\hat{\tb{x}}_n^C(\tb{x}_0){\vphantom{\hat{\tb{W}}}}^{(\boldsymbol{r}-\boldsymbol{r'})}\mathds{E}\big[\Delta\hat{\tb{W}}^{\boldsymbol{r'}}_{n,N}(\tb{x}_0)\big],
\end{equation}
for any multi-index $\tb{r}$. Because of the truncation of order $N$, the approximation is accurate for moments of order $k\leq N$. For example, if we choose $N=1$, the approximation is truncated to order $1$ and the polynomial is linear in the noise and hence, not dependent on the second or higher moments of the r.p. $\tb{W}_n$.}

\cblack{The convergence of the polynomial expansions presented above is rigorously established in Section \ref{var_sol_noise}.}


\subsection{Initial uncertainty}\label{PCE_subsec}

\cblack{
In general, the initial condition for the SDE \eqref{dyn_sde} is unknown and $\tb{X}_0$ is modelled as a random vector with a given probability distribution. It is tempting to handle this uncertainty as an initial effective noise, i.e., to assume that $\Delta\hat{\tb{W}}_{0,N}(\tb{x}_0)=\tb{X}_0$ and $\tb{x}_0=\textbf{0}$ in Eqs. \eqref{eff_noise_expansion} and \eqref{cond_moms_eff_noise}. However, this approach turns out naive. Since Eqs. \eqref{eff_noise_expansion} and \eqref{cond_moms_eff_noise} are obtained from a Taylor expansion of $\Delta\hat{\tb{W}}_{0,N}(\tb{x}_0)$, when the higher order moments of the effective noise are significant we need to increase the order $N$ of the approximation in order to maintain a prescribed (sufficiently good) accuracy. A larger $N$ implies the computation of higher-order derivatives of functions $\tb{u}$ and $\tb{G}$ and, as a consequence, an increased computational effort. In general, the uncertainty of the initial conditions can be expected to be independent of the dynamical perturbation $\tb{W}(t)$ and, possibly, to have a larger power and more significant higher-order moments compared to the process $\tb{W}(t)$. For these reasons, it is more convenient to handle the initial uncertainty using a specific expansion of order possibly higher than $N$.
}

\cblack{
The polynomial chaos expansion (PCE) method \cite{Lu16} is a technique that provides a polynomial expansion of a r.v. propagated through a deterministic dynamical system. The standard PCE scheme cannot be used in a SDE like Eq. \eqref{dyn_sde}. However, the argument in Section \ref{pol_eff_noise_A} enables us to circumvent this problem, as we have already obtained a deterministic recursion for the moments of the effective noise in Eq. \eqref{cond_moms_eff_noise}.
}

\cblack{
In order to compute a PCE of the conditional moments of $\hat{\tb{X}}_n$ we take a set of $N_p$ polynomials $\{\Phi_i :\mathds{R}^v\rightarrow\mathds{R} \}_{i=1}^{N_p}$, selected to be orthogonal with respect to the pdf $f_{\tb{X}_0}$ of the initial condition $\tb{X}_0$. Then, we construct the approximation
\begin{equation}\label{pol_exp_1}
\mathds{E}\big[\hat{\tb{X}}^{\tb{r}}_{n,N}\big|\mX_0\big]\simeq\sum_{i=1}^{N_p}c^{(\tb{r})}_{i,n,N}\Phi_i(\tb{X}_0)\,,
\end{equation}
where $\mX_0$ is the $\sigma$-algebra generated by $\tb{X}_0$ and the $c^{(\tb{r})}_{i,n,N}$'s are the PCE coefficients (note thet the superscript ${}^{(\tb{r})}$ simply indicates dependence on the multi-index $\tb{r}$ on the left-hand side). A simple way to compute these coefficients is the so-called non-intrusive method \cite{Lu16}, for which
\begin{equation}\label{minim_PCE}
\left\{c^{(\tb{r})}_{i,n,N}\right\}_{i=1}^{N_p}=\underset{\;\{c_k\}_{k=1}^{N_p}}{\textrm{argmin}}\int\left(\mathds{E}\big[\hat{\tb{X}}^{\boldsymbol{r}}_{n,N}(\tb{x}_0)\big]-\sum_{k=1}^{N_p} c_k\Phi_k(\tb{x}_0)\right)^{\!2}\!\!f_{\tb{X}_0}(\tb{x}_0)\diff\tb{x}_0\,.
\end{equation}
While the optimisation problem \eqref{minim_PCE} above cannot be solved exactly in general, for most practical applications it is possible to approximate the integral using Monte Carlo. If we draw $N_s$ samples from the pdf $f_{\tb{X}_0}$, denoted by $\tb{X}_{0,j}$, $j=1,\ldots,N_s$, it is straightforward to compute an approximation of the PCE coefficients by solving the linear least-squares problem
\begin{equation}\label{minim_PCE_disc}
\left\{\hat{c}^{(\tb{r})}_{i,n,N}\right\}_{i=1}^{N_p}=\underset{\;\{c_k\}_{k=1}^{N_p}}{\textrm{argmin}}\sum_{j=1}^{N_s}\left(\mathds{E}\big[\hat{\tb{X}}^{\boldsymbol{r}}_{n,N}(\tb{X}_{0,j})\big]-\sum_{k=1}^{N_p} c_k\Phi_k(\tb{X}_{0,j})\right)^{\!2}\!\!,
\end{equation}
which, in turn, yields the approximate conditional moments
\begin{equation}\label{momet_approx_PCE}
\mathds{E}\big[\hat{\tb{X}}^{\tb{r}}_{n,N}\big|\mX_0\big]\simeq\mathds{E}\big[\hat{\tb{X}}^{\tb{r}}_{n,N}\big|\mX_0\big]_{N_p}\!\!\coloneq\sum_{i=1}^{N_p}\hat{c}^{(\tb{r})}_{i,n,N}\Phi_i(\tb{X}_0)\,.
\end{equation}
}

\cblack{
Some remarks are in order regarding the validity of Eq.\,\eqref{momet_approx_PCE}:
}

\cblack{
\begin{itemize}
\item When the dimension of the state space is $v$, the PCE approximation is of order $N_{\textrm{PCE}}$ where the number of orthogonal polynomials is \cite[Chapter 2]{Be99}
\begin{equation}\label{N_p_coeffs}
N_p=\begin{pmatrix}
N_{\textrm{PCE}}+v\\
v
\end{pmatrix}.
\end{equation}
\item The least-squares problem in \eqref{minim_PCE_disc} can be solved when the correlation matrix $\boldsymbol{\Phi}^\top\boldsymbol{\Phi}$ has full rank, where
\begin{equation}\label{matrix_ortho}
\boldsymbol{\Phi}\coloneq\begin{pmatrix}
					\Phi_1(\tb{X}_{0,1})&\cdots&\Phi_{N_p}(\tb{X}_{0,1})\\
					\vdots&\ddots&\vdots\\
					\Phi_1(\tb{X}_{0,N_s})&\cdots&\Phi_{N_p}(\tb{X}_{0,N_s})
				\end{pmatrix}.
\end{equation}
This implies that $N_s\geq N_p$ (in practice, $N_s>N_p$ and sufficiently large). The numerical computation of \eqref{minim_PCE_disc} is typically more stable when the polynomials $\{\Phi_i\}_{i=1}^{N_p}$ are orthonormal \cite{DE87}, i.e., when $\|\Phi_i\|=1$.
\item The polynomial expansions \eqref{pol_exp_1} and \eqref{momet_approx_PCE} converge in mean square error (MSE) when the $L^2$-norm of the conditional moments with respect to the density $f_{\tb{X}_0}(\tb{x}_0)\diff\tb{x}_0$ are finite \cite{Ch11}, i.e.,
\begin{equation}\label{L2_condition}
\int \mathds{E}\big[\hat{\tb{X}}^{\boldsymbol{r}}_{n,N}(\tb{x}_0)\big]^2f_{\tb{X}_0}(\tb{x}_0)\diff\tb{x}_0<\infty.
\end{equation}
for the selected multi-index $\tb{r}$.
\end{itemize}
}

\cblack{
Finally, we recall the rule of iterated expectations \cite[Theorem 3.4]{Wa04}), which yields
\begin{equation}\label{expectation_sol_cond_total}
\mathds{E}\big[\hat{\tb{X}}^{\boldsymbol{r}}_{n,N}\big]=\mathds{E}\Big[\mathds{E}\big[\hat{\tb{X}}^{\boldsymbol{r}}_{n,N}\big|\mX_0\big]\Big]\simeq\sum_{i=1}^{N_p}\hat{c}^{(\tb{r})}_{i,n,N}\mathds{E}\big[\Phi_i(\tb{X}_0)\big]\,.
\end{equation}
When the polynomials $\Phi_i$ are orthonormal with respect to $f_{\tb{X}_0}$ it follows that
\begin{equation}
\int\limits_{\;\;\;\mathds{R}}\Phi_i(\tb{x}_0)\Phi_j(\tb{x}_0)f_{\tb{X}_0}(\tb{x}_0)\diff\tb{x}_0=\delta_{ij} := \left\{
	\begin{array}{ll}
	1, &\text{if $i=j$},\\
	0, &\text{otherwise},\\
	\end{array}
\right.
\end{equation}
and in particular,
\begin{equation}
\mathds{E}\big[\Phi_i(\tb{X}_0)\big]=\delta_{i1}.
\end{equation}
Therefore, Eq.\,\eqref{expectation_sol_cond_total} readily yields
\begin{equation}\label{moments_PCE}
\mathds{E}\big[\hat{\tb{X}}^{\tb{r}}_{n,N}\big]\simeq\mathds{E}\big[\hat{\tb{X}}^{\tb{r}}_{n,N}\big]_{N_p}=\hat{c}_{1,n,N}^{(\tb{r})}
\end{equation}
when the expansion is based on an orthonormal set of polynomials.
}



\subsection{1-Dimensional marginal densities} \label{Marg_ps}

\cblack{
The approximate moments in Eq.\,\eqref{moments_PCE} yield a partial description of the probability distribution of $\hat{\tb{X}}_n$. However, in many problems, the uncertainty associated to the random sequence $\hat{\tb{X}}_n$ is easier to interpret in terms of the probability density function (pdf) of the r.v.'s of interest. In this section, we describe a procedure to approximate the marginal pdf of each variable $\hat{X}^{(k)}_n$ using a Gram--Charlier expansion \cite{Ko97}.
}

\cblack{
We introduce some notation first. Let $X$ be a real r.v. The pdf of $X$ is denoted by $f_{X}$, while $\Psi_{X}(t)\coloneq\mathds{E}[\e^{iXt}]$ is the characteristic function of $X$ (where $i$ is the imaginary unit and $t \in \mathds{R}$). The cumulant generating function of $X$ is \cite{Ko97}
\begin{equation}\label{def_log_char}
\log\big(\Psi_{X}(t)\big)=\sum_{r=1}^\infty\frac{\kappa_{r}(X)}{r!}(it)^{r},
\end{equation}
where $\kappa_{r}$ denotes the $r$-th order cumulant. The cumulants $\kappa_{r}(X)$ can be computed in terms of the moments of $X$ using lookup tables \cite{Ko97}.
}

\cblack{
The general Gram--Charlier expansion of the marginal pdf of $\hat{X}^{(k)}_{n,N}$ conditional on a fixed initialization $\tb{X}_0=\tb{x}_0$ can be written as
\begin{equation}\label{char_gen_form}
f_{\hat{X}^{(k)}_{n,N}|\tb{x}_0}(x|\tb{x}_0)\simeq\left[1+\sum_{r=1}^{N}\frac{(-1)^{r}}{r!}C_{r}\big[\hat{X}^{(k)}_{n,N}(\tb{x}_0),Z_\varphi\big]\frac{\diff^{r}}{\diff {x}^{r}}\right]\!\varphi(x)
\end{equation}
for any $x\in\mathds{R}$ (such that the expansion converges), where 
\begin{multline}\label{Coeff_GC}
C_{r}\big[\hat{X}^{(k)}_{n,N}(\tb{x}_0),Z_\varphi\big] \coloneq 
B_{r}\Big(\kappa_1\big(\hat{X}^{(k)}_{n,N}(\tb{x}_0)\big)-\kappa_1(Z_\varphi), \ldots\\
\ldots,\kappa_r\big(\hat{X}^{(k)}_{n,N}(\tb{x}_0)\big)-\kappa_r(Z_\varphi)\Big),
\end{multline}
and $B_{r}$ is the $r$-th Bell polynomial \cite{Ab05}, $\varphi$ is an auxiliary pdf and $Z_\varphi$ is a r.v. with density $\varphi$.
}

\cblack{
In the case at hand, we note that for a fixed $\tb{X}_0=\tb{x}_0$ the distribution of the solution $\hat{\tb{X}}_n(\tb{x}_0)=\hat{\tb{x}}_n^C(\tb{x}_0)+\Delta\hat{\tb{W}}_n(\tb{x}_0)$ depends essentially on the distribution of the effective noise $\Delta\hat{\tb{W}}_n(\tb{x}_0)$, as $\hat{\tb{x}}_n^C(\tb{x}_0)$ is the numerical approximation of the deterministic solution to the ODE $\dot{\tb{x}}^{C}(t)=\boldsymbol{u}(\tb{x}^{C},t)$ with initial condition $\tb{x}^{C}(t_0)=\tb{x}_0$. Using Eq.\,\eqref{eff_noise_expansion}, we can expand the effective noise in terms of the noise increments $\Delta\tb{W}_{n}$. Specifically, if we apply a truncation of order $N=1$, the $k$-th effective noise coordinate becomes
\begin{equation}\label{expansion_N_1}
\Delta\hat{W}_{n,1}^{(k)}(\tb{x}_0)=\sum_{j=1}^d\sum_{m=1}^n b_{n,m,j}^{(k)}(\tb{x}_0)\Delta W_{m}^{(j)},
\end{equation}
where the $b_{n,m,j}^{(k)}(\tb{x}_0)$'s are deterministic coefficients. Hence, $\Delta\hat{W}_{n,1}^{(k)}(\tb{x}_0)$ is a linear combination of independent r.v.'s. If $\tb{W}(t)$ is a Wiener process, then $\Delta\hat{W}_{n,1}^{(k)}(\tb{x}_0)$ is Gaussian and, even for more general processes, recent results on Berry--Esseen bounds \cite{La10,Kl10} suggest that a Gaussian approximation for $\Delta\hat{W}_{n,1}^{(k)}(\tb{x}_0)$ is a plausible choice. Therefore, we let the auxiliary pdf $\varphi(x)$ in Eq.\,\eqref{char_gen_form} be a normal pdf depending on $\tb{x}_0$, denoted by $\varphi^{G^{(k)}}_{n}(x|\tb{x}_0)$ with mean $\mu^{(k)}_{n}(\tb{x}_0)$ and standard deviation $\sigma^{(k)}_{n}(\tb{x}_0)$. We write $Z_{\varphi}^{G^{(k)}}\!(\tb{x}_0)$ to denote a r.v. with pdf precisely $\varphi^{G^{(k)}}_{n}(\cdot\,|\tb{x}_0)$.
}

\cblack{
The Gram--Charlier expansion with a Gaussian auxiliary density is well studied and known as Gram--Charlier expansion of type $A$. In particular, Eq. \eqref{char_gen_form} can be rewritten as \cite{Ko97}
\begin{multline}\label{GC_Type_A}
f_{\hat{X}^{(k)}_{n,N}|\tb{x}_0}(x|\tb{x}_0)\simeq\left[1+\sum_{r=1}^N\frac{1}{r!\,\sigma^{(k)}_{n}(\tb{x}_0)^{\vphantom{A}^{r}}}C_{r}\Big[\hat{X}^{(k)}_{n,N}(\tb{x}_0),Z_{\varphi}^{G^{(k)}}\!(\tb{x}_0)\Big]\vphantom{H_{e_r}\left(\frac{x-\mu^{(k)}_{n}(\tb{x}_0)}{\sigma^{(k)}_{n}(\tb{x}_0)}\right)}\right.
\\\left.\vphantom{1+\sum_{r=1}^N\frac{1}{r!(\sigma^{(k)}_{n}(\tb{x}_0))^{\vphantom{A}^{r}}}C_{r}\big[\hat{X}^{(k)}_{n,N}(\tb{x}_0)\,;\varphi^G_{n,\tb{x}_0}\big]} \times H_{e_r}\left(\frac{x-\mu^{(k)}_{n}(\tb{x}_0)}{\sigma^{(k)}_{n}(\tb{x}_0)}\right)\right]\varphi^{G^{(k)}}_{n}\!(x|\tb{x}_0),
\end{multline}
where $H_{e_r}(x)$ is the $r$-th Hermite polynomial that satisfies the Rodrigues formula \cite{Ra81}
\begin{equation}\label{Rodrigues}
H_{e_r}(x) =(-1)^{r}\e^{x^2/2}\frac{\diff^{r}}{\diff x^{r}}\e^{-x^2/2}.
\end{equation}
}

\cblack{
For simplicity, we propose to compute the mean $\mu^{(k)}_{n}(\tb{x}_0)$ and standard deviation $\sigma^{(k)}_{n}(\tb{x}_0)$ of the auxiliary Gaussian density $\varphi^{G^{(k)}}_{n}\!(x|\tb{x}_0)$ using the truncations of order $N=1$ of the effective noise $\Delta\hat{W}_{n,1}^{(k)}(\tb{x}_0)$. This yields
\begin{equation}\label{mu_gauss}
\mu^{(k)}_{n}(\tb{x}_0)=\hat{\tb{x}}_{n}^{C^{(k)}}\!(\tb{x}_0)+\mathds{E}\big[\Delta\hat{W}_{n,1}^{(k)}(\tb{x}_0)\big],
\end{equation}
and
\begin{equation}\label{sigma_gauss}
\sigma^{(k)}_{n}(\tb{x}_0)=\sqrt{\mathds{E}\big[\big(\Delta\hat{W}_{n,1}^{(k)}(\tb{x}_0)\big)^{2}\big]-\mathds{E}\big[\Delta\hat{W}_{n,1}^{(k)}\big(\tb{x}_0)\big]^2},
\end{equation}
for $k=1,\ldots,v$ where $\mathds{E}\big[\Delta\hat{W}_{n,1}^{(k)}\big(\tb{x}_0)\big]$ and $\mathds{E}\big[\big(\Delta\hat{W}_{n,1}^{(k)}(\tb{x}_0)\big)^{2}\big]$ can be computed recursively from Eq.\,\eqref{eff_noise_expansion}.
}

\cblack{
The convergence of the expansion in Eq.\,\eqref{GC_Type_A}, i.e., the approximation error when the series is truncated to some finite order is addressed in Section \ref{pdf_sol_noise}.
}

\cblack{
When the initial condition $\tb{X}_0$ is random, it is possible to construct PCE approximations of $\mu^{(k)}_{n}$, $\sigma^{(k)}_{n}$ and $C_{r}$ in a similar way as we computed the approximations of conditional moments of $\hat{\tb{X}}^{\boldsymbol{r}}_{n,N}$ in Section \ref{PCE_subsec}. In particular,
\begin{equation}\label{mu_sigma_PCE}
\mu^{(k)}_{n}(\tb{X}_0)_{N_p}\!\!=\sum_{i=1}^{N_p}\hat{c}_{i,n}^{\mu(k)}\Phi_{i}\big(\tb{X}_0\big),\qquad\sigma^{(k)}_{n}(\tb{X}_0)_{N_p}\!\!=\sum_{i=1}^{N_p}\hat{c}_{i,n}^{\sigma(k)}\Phi_{i}\big(\tb{X}_0\big),
\end{equation}
and
\begin{equation}\label{Bell_PCE}
C_{r}\Big[\hat{X}^{(k)}_{n,N}(\tb{X}_0),Z_\varphi^{G^{(k)}}\!(\tb{X}_0)\Big]_{N_p}\!\!=\sum_{i=1}^{N_p}\hat{c}_{i,n}^{C_{r}(k)}\Phi_{i}\big(\tb{X}_0\big),
\end{equation}
where, the same as in Section \ref{PCE_subsec}, the coefficients of the expansion are obtained by solving the least-squares problem
\begin{equation}\label{minim_PCE_for_pdf}
\left\{\hat{c}^{\,[s](k)}_{i,n}\right\}_{i=1}^{N_p}=\underset{\;\{c_k\}_{k=1}^{N_p}}{\textrm{argmin}}\sum_{j=1}^{N_s}\left(u^{[s](k)}_{n}(\tb{X}_{0,j})-\sum_{k=1}^{N_p} c_k\Phi_k(\tb{X}_{0,j})\right)^{\!2}\!\!,\qquad s=1,2,3,
\end{equation}
where 
\begin{itemize}
\item $\hat{c}^{\,[1](k)}_{i,n}=\hat{c}_{i,n}^{\mu(k)}$ and $u^{[1](k)}_{n}=\mu^{(k)}_{n}$;
\item $\hat{c}^{\,[2](k)}_{i,n}=\hat{c}_{i,n}^{\sigma(k)}$ and $u^{[2](k)}_{n}=\sigma^{(k)}_{n}$;
\item $\hat{c}^{\,[3](k)}_{i,n}=\hat{c}_{i,n}^{\,C_{r}}$ and $u^{[3](k)}_{n}=C_{r}\Big[\hat{X}^{(k)}_{n,N},Z_\varphi^{G^{(k)}}\!\Big]$.
\end{itemize}
}

\cblack{
Finally, if we draw $N'_s$ i.i.d. samples from the random initial condition $\tb{X}_0$, denoted $\tb{X}'_{0,j}$, $j=1,\ldots,N'_s$, then we can use Eq. \eqref{GC_Type_A} to approximate the pdf of $\hat{X}^{(k)}_{n,N}$ as
\begin{equation}\label{law_tot_prob}
f_{\hat{X}^{(k)}_{n,N}}(x)=\!\!\!\!\int\limits_{\quad\mathds{R}^v}\!\!f_{\hat{X}^{(k)}_{n,N}|\tb{X}_0}(x|\tb{x}_0)f_{\tb{X}_0}(\tb{x}_0)\diff\tb{x}_0\simeq\frac{1}{N'_s}\sum_{j=1}^{N'_s}f_{\hat{X}^{(k)}_{n,N}|\tb{X}_0}(x|\tb{X}'_{0,j}).
\end{equation}
}


\subsection{Outline of the algorithms}\label{sec_pseudocode}

\cblack{
In this section we provide a summary of the proposed algorithms for the approximate computation of the moments $\mathds{E}\big[\hat{\tb{X}}^{\boldsymbol{r}}_{n,N}\big]$ and the 1-dimensional marginal densities $f_{\hat{X}^{(k)}_{n,N}}(x)$, $k=1, \ldots, v$.
}

\cblack{
Table \ref{MAIN_proc} provides a list, with brief descriptions, of the inputs and outputs of the two proposed approximation schemes. Algorithm \ref{alg_main_ps_1} displays a pseudocode, with cross-references to Sections \ref{pol_eff_noise_A} and \ref{PCE_subsec}, of the numerical scheme for the computation of moments assuming a random initial condition $\tb{X}_0$. If $\tb{X}_0=\tb{x}_0$ the algorithm is simply run with $N_s=1$. Algorithm \ref{alg_main_ps_2} shows a pseudocode for the approximation of marginal densities, with cross-references to Section \ref{Marg_ps}. 
}

\begin{table}[htb]
\begin{center}
	\vspace{0.cm}\begin{tabular}{|m{2cm}||m{8.7cm}|}
	\hline
	\hfil\textbf{Inputs} & \hfil\textbf{Description}\\	
	\hline\hline
	\hfil$h$ &Step-size.\\
	\hline
	\hfil$t_0$ &Initial time.\\
	\hline
	\hfil$t_n$ &Final time.\\
	\hline
	\hfil$N$ &Order of the polynomial expansions.\\
	\hline
	\hfil$v$ &Dimension of $\tb{X}(t)$.\\
	\hline
	\hfil$d$ &Dimension of $\tb{W}(t)$.\\
	\hline
	\hfil$\tb{X}_0$ &Initial condition.\\
	\hline
	\hfil$\vphantom{A^{B^B}_{\displaystyle{A}}}\;\mathds{E}\big[\Delta\tb{W}_{m}^{\,\boldsymbol{r}_{\vphantom{a_1}}}\big]\;$ &Moments of the noise increments for $m\geq 1$ and $|\boldsymbol{r}|\leq N$.\\
	\hline
	\hfil$N_{\textrm{PCE}}$ &Truncation order of the PCE scheme (when $\tb{X}_0=\tb{x}_0$ is fixed, this is not needed).\\
	\hline
	\hfil$N_s$ &Number of samples (if $\tb{X}_0=\tb{x}_0$ is fixed, $N_s=1$).\\
	\hline
	\hfil$N'_s$ &Number of i.i.d. samples of $\tb{X}_0$ to approximate the pdf $f_{\hat{X}^{(k)}_{n,N_{\vphantom{B_A}}}}(x)$ in Eq. \eqref{law_tot_prob}.\\
	\hline
	\hfil$\tb{u}$ &\textit{Drift coefficient} in Eq.\,\eqref{dyn_sde}.\\
	\hline
	\hfil$\tb{G}$ &\textit{Diffusion coefficient} in Eq.\,\eqref{dyn_sde}.\\
	 \hline
	 \multicolumn{1}{c}{}&\multicolumn{1}{c}{} \vspace{-0.2cm}\\
	 \hline
	\hfil\textbf{Outputs} & \hfil\textbf{Description}\\
	\hline\hline
	\hfil$\;\mathds{E}\big[\hat{\tb{X}}^{\boldsymbol{r}}_{n,N}\big]_{N_p}$ &Moments of the numerical solution of Eq. \eqref{dyn_sde}, for $|\tb{r}|\leq N$, computed with a basis of $N_p$ orthogonal polynomials, where $N_p$ is given by Eq. \eqref{N_p_coeffs}.\vspace{0.02cm}\\
	\hline
	\hfil$f_{\hat{X}^{(k)}_{n,N}}$ &
	Estimate of the pdf $f_{\hat{X}_n}$ where $\hat{\tb{X}}_n$ is the numerical approximation of the r.v. $\tb{X}(t_n)$\\
\hline	
\end{tabular}
\end{center}
\caption{Inputs and outputs of the algorithms for moment computation and estimation of 1-dimensional marginal pdf's.}
\label{MAIN_proc}
\end{table}

\vspace{0.cm}
\begin{algorithm}
\caption{Computation of moments}
\label{alg_main_ps_1}
\begin{algorithmic}[1]
\STATE Generate $N_s$ samples of $\tb{X}_0$, denoted $\tb{X}_{0,j}$, $j=1, \ldots, N_s$.
\STATE Compute $N_p$ using Eq. \eqref{N_p_coeffs} and matrix $\boldsymbol{\Phi}$ using Eq. \eqref{matrix_ortho}\big) such that $\boldsymbol{\Phi}^\top\boldsymbol{\Phi}$. We assume $\boldsymbol{\Phi}$ is full-rank.
\STATE Set $\hat{\tb{x}}_0^C(\tb{X}_{0,j})=\tb{X}_{0,j}$, $\Delta\hat{\tb{W}}_{0,N}(\tb{X}_{0,j})=0$ and $\Delta\hat{\tb{W}}_{0,1}(\tb{X}_{0,j})=0$ for $j=1,\ldots,N_s$.
\STATE Set $n=\lceil(t_n-t_0)/h\rceil$, where $\lceil\,\cdot\,\rceil$ denotes the ceiling function.
\FOR{$m=1,\ldots,n$}
	\FOR{$j=1,\ldots,N_s$}
		\STATE Evaluate the central part $\hat{\tb{x}}_m^C(\tb{X}_{0,j}) = \hat{\tb{x}}_{m-1}^C(\tb{X}_{0,j}) + h\boldsymbol{u}\big(\hat{\tb{x}}_{m-1}^C(\tb{X}_{0,j}),t_{m-1}\big),$ where $t_m=t_0+mh$.
		\phantomsection\label{line_8}\STATE Evaluate $\mathds{E}\big[\Delta\hat{\tb{W}}^{\tb{r}}_{m,N}(\tb{X}_{0,j})\big]$ for $1\leq|\boldsymbol{r}|\leq N$ using Eq.\,\eqref{cond_moms_eff_noise}.
		\phantomsection\label{line_9}\STATE Evaluate $\mathds{E}\big[\Delta\hat{W}_{m,1}^{(k)}(\tb{X}_{0,j})\big]$ and $\mathds{E}\big[\big(\Delta\hat{W}_{m,1}^{\,\scriptscriptstyle{(k)}_{\vphantom{a_1}}}(\tb{X}_{0,j})\big)^{2}\big]$ for $k=1,\ldots,v$.
	\ENDFOR
\ENDFOR
\STATE Compute $\mathds{E}\big[\hat{\tb{X}}^{\boldsymbol{r}}_{n,N}\left(\tb{X}_{0,j}\right)\big]$ for  $j=1,\ldots,N_s$ using Eq.\,\eqref{cond_moms_sol}.
\STATE Solve the least-squares problem \eqref{minim_PCE_disc} and compute $\mathds{E}\big[\hat{\tb{X}}^{\boldsymbol{r}}_{n,N}\big]_{N_p}$ for $1\leq|\boldsymbol{r}|\leq N$ with Eq.\,\eqref{moments_PCE}.
\end{algorithmic}
\end{algorithm}
\vspace{1cm}

\begin{algorithm}
\caption{Computation of 1-dimensional marginal pdf's}
\label{alg_main_ps_2}
\begin{algorithmic}[1]
\STATE Generate $N_s$ samples of $\tb{X}_0$, denoted $\tb{X}_{0,j}$, $j=1, \ldots, N_s$.
\FOR{$j=1,\ldots,N_s$}
\FOR{$k=1,\ldots,v$}
\STATE Compute $\mu^{\scriptscriptstyle{(k)}}_{n}(\tb{X}_{0,j})$ and $\sigma^{\scriptscriptstyle{(k)}}_{n}(\tb{X}_{0,j})$ using Eqs.\,\eqref{mu_gauss} and \eqref{sigma_gauss} respectively.
\STATE Compute $C_{r}\Big[\hat{X}^{\scriptscriptstyle{(k)}}_{n,N}(\tb{X}_{0,j}),Z_\varphi^{G^{(k)}}\!(\tb{X}_{0,j})\Big]$ for $r=1,\ldots,N$.
\ENDFOR
\ENDFOR
\STATE Solve the least-squares problem to compute the PCE coefficients of Eqs.\,\eqref{mu_sigma_PCE} and \eqref{Bell_PCE}.
\STATE Generate $N'_s$ samples of $\tb{X}_0$, denoted $\tb{X}'_{0,j}$, $j=1, \ldots, N'_s$.
\FOR{$j=1,\ldots,N'_s$}
\FOR{$k=1,\ldots,v$}
\STATE Compute $\mu^{\scriptscriptstyle{(k)}}_{n}(\tb{X}'_{0,j})$ and $\sigma^{\scriptscriptstyle{(k)}}_{n}(\tb{X}'_{0,j})$ using Eq.\,\eqref{mu_sigma_PCE}.
\STATE Compute $C_{r}\Big[\hat{X}^{\scriptscriptstyle{(k)}}_{n,N}(\tb{X}'_{0,j}),Z_{\varphi}^{G^{(k)}}\!(\tb{X}'_{0,j})\Big]$, for $r=1,\ldots,N$, using Eq.\,\eqref{Bell_PCE}.
\ENDFOR
\ENDFOR
\STATE Compute the coefficients of the Hermite polynomials $H_{e_r}$ for $r=0,\ldots,N$.
\STATE Apply Eq.\,\eqref{law_tot_prob}, combined with Eq.\,\eqref{GC_Type_A}, to compute $f_{\hat{X}^{\scriptscriptstyle{(k)}}_{n,N}}(x)$ for any $x\in\mathds{R}$ and $k=1,\ldots,v$.
\end{algorithmic}
\end{algorithm}
\vspace{0.2cm}


\section{Variational solution based in polynomial expansion over the noise}\label{var_sol_noise}

\cblack{
In this section we provide the analysis needed to support the results in Section \ref{pol_eff_noise_A} and, specifically, Algorithm \ref{alg_main_ps_1} for the approximate computation of moment of the random sequence $\hat X_n$. Our analysis relies on the notion of convergence region for a Taylor expansion as defined below.
}

\cblack{
\begin{definition}\label{Conv_reg_def}
Let $\boldsymbol{\mathrm{g}}$ be a smooth function, $\boldsymbol{\mathrm{g}}:\mathds{R}^v\times\mathds{R}_0^{+}\rightarrow\mathds{R}$. The convergence region of the Taylor expansion of $\mathrm{g}$, centred around $\tb{x}_0\in\mathds{R}^v$ at time $t\in\mathds{R}_0^{+}$, is the set 
$$
\rho_{\mathrm{g}}(\tb{x}_0,t)\coloneq\Big\{\tb{x}\in\mathds{R}^v:\sum_{|\boldsymbol{r}|=0}^\infty\frac{1}{\boldsymbol{r}!}\frac{\partial^{|\boldsymbol{r}|}}{\partial \tb{x}^{\boldsymbol{r}}}\mathrm{g}(\tb{x}_0,t)\tb{x}^{\boldsymbol{r}}<\infty\Big\} \subseteq \mathds{R}^v.
$$
\end{definition}
}


\cblack{Let us assume a fixed initial condition $\tb{X}_0=\tb{x}_0$. The moments of the sequence $\hat X_n$ follow readily from the statistics of the effective noise sequence $\Delta\hat{W}_{n}^{\,\scriptscriptstyle{(k)}_{\vphantom{a_1}}}(\tb{x}_0)$. Therefore, we start with the expansion formula for the effective noise in Eq. \eqref{eff_noise_expansion}.
}

\cblack{
\begin{theorem}\label{eff_noise_approx}
Assume that the functions $\boldsymbol{u}$ and $\boldsymbol{G}$ in Eq. \eqref{dyn_sde} are real and smooth. For any positive integers $N$ and $n$, the effective noise given an initial condition $\tb{X}_0=\tb{x}_0$ can be written as:
\begin{eqnarray}
\Delta\hat{W}_{n}^{\,\scriptscriptstyle{(k)}_{\vphantom{a_1}}}(\tb{x}_0) &=& \Delta\hat{W}_{n-1}^{\,\scriptscriptstyle{(k)}_{\vphantom{a_1}}}(\tb{x}_0)
+h\!\!\sum_{|\boldsymbol{r}|=1}^N\frac{1}{\boldsymbol{r}!}\frac{\partial^{|\boldsymbol{r}|}}{\partial \tb{x}_{n-1}^{\boldsymbol{r}}}u^{\scriptscriptstyle{(k)}_{\vphantom{a_1}}}\big(\hat{\tb{x}}_{n-1}^C(\tb{x}_0),t_{n-1}\big)\Delta\hat{\tb{W}}_{n-1}^{\boldsymbol{r}}(\tb{x}_0) \nn \\
&& + \sum_{|\boldsymbol{r}|=0}^{N-1}\sum_{j=1}^d\frac{1}{\boldsymbol{r}!}\frac{\partial^{|\boldsymbol{r}|}}{\partial \tb{x}_{n-1}^{\boldsymbol{r}}}G^{(k,j)}\big(\hat{\tb{x}}_{n-1}^C(\tb{x}_0),t_{n-1}\big)\Delta W^{(j)}_{n}\Delta\hat{\tb{W}}_{n-1}^{\boldsymbol{r}}(\tb{x}_0) \nn \\
&& + R^{(k)}_{n,N}\big(\Delta\tb{W}_n,\Delta\hat{\tb{W}}_{n-1}(\tb{x}_0)\big), \label{eqJ-1} 
\end{eqnarray}
where $\Delta\hat{\tb{W}}_{0}(\tb{x}_0)=0$ and $R^{\scriptscriptstyle{(k)}}_{n,N}$ is the remainder term of the polynomial expansion at step $n$ with truncation order $N$. If 
\begin{equation}
\Delta\hat{\tb{W}}_{n-1}(\tb{x}_0) \in
\rho_{u^{(k)}}\big(\hat{\tb{x}}^C_{n-1}(\tb{x}_0),t_{n-1}\big)
\cap_{j=1}^d 
\rho_{G^{(k,j)}}\big(\hat{\tb{x}}^C_{n-1}(\tb{x}_0),t_{n-1}\big), \nn
\end{equation} 
then
$$
\lim_{N\rightarrow\infty}R^{(k)}_{n,N}\big(\Delta\tb{W}_n,\Delta\hat{\tb{W}}_{n-1}(\tb{x}_0)\big)=0.
$$  
\end{theorem}
}

\cblack{
\begin{remark}
Note that $\frac{\Delta \hat{\tb{W}}_{n-1}^{\tb{r}}(\tb{x}_0)}{\tb{r}!} = 1$ for $|\tb{r}|=0$.
\end{remark}
}

\cblack{
{\sc Proof:} Recall the decomposition of the sequence $\hat{\tb{X}}_n(\tb{x}_0)$ into its central part and the effective noise,
\begin{equation}
\hat{\tb{X}}_{n-1}(\tb{x}_0)=\hat{\tb{x}}_{n-1}^C(\tb{x}_0)+\Delta\hat{\tb{W}}_{n-1}(\tb{x}_0).
\label{eqJ0}
\end{equation}
Using the relationship above, the the Taylor expansions of $u^{(k)}$ (of order $N$) and $G^{(kj)}$ (of order $N-1$) with respect to $\hat{\tb{X}}_{n-1}(\tb{x}_0)$ and centred at $\hat{\tb{x}}_{n-1}^C(\tb{x}_0)$ can be written as
\begin{eqnarray}
u^{(k)}\big(\hat{\tb{X}}_{n-1},t_{n-1}\big)
&=&
\sum_{|\boldsymbol{r}|=0}^N\frac{1}{\boldsymbol{r}!}\frac{\partial^{|\boldsymbol{r}|}}{\partial \tb{x}_{n-1}^{\boldsymbol{r}}}u^{(k)}\big(\hat{\tb{x}}_{n-1}^C(\tb{x}_0),t_{n-1}\big)\Delta\hat{\tb{W}}_{n-1}^{\boldsymbol{r}}(\tb{x}_0) \nn\\
&&+R^{u^{(k)}}_{n-1,N}\big(\Delta\hat{\tb{W}}_{n-1}(\tb{x}_0)\big) 
\label{eqJ1}
\end{eqnarray}
and
\begin{eqnarray}
G^{(k,j)}\big(\hat{\tb{X}}_{n-1},t_{n-1}\big) &=& \sum_{|\boldsymbol{r}|=0}^{N-1}\frac{1}{\boldsymbol{r}!}\frac{\partial^{|\boldsymbol{r}|}}{\partial \tb{x}_{n-1}^{\boldsymbol{r}}}G^{(k,j)}\big(\hat{\tb{x}}_{n-1}^C(\tb{x}_0),t_{n-1}\big)\Delta\hat{\tb{W}}_{n-1}^{\boldsymbol{r}}(\tb{x}_0) \nn\\
&&+R^{G^{(k,j)}}_{n-1,N-1}\big(\Delta\hat{\tb{W}}_{n-1}(\tb{x}_0)\big),
\label{eqJ2}
\end{eqnarray}
respectively, where $R^{u^{(k)}}_{n-1,N}\big(\Delta\hat{\tb{W}}_{n-1}(\tb{x}_0)\big)$ and $R^{G^{(k,j)}}_{n-1,N-1}\big(\Delta\hat{\tb{W}}_{n-1}(\tb{x}_0)\big)$ are remainder terms. If we substitute Eqs. \eqref{eqJ0}--\eqref{eqJ2} into Euler--Maruyama scheme of \eqref{sto_integrator}, we obtain the expansion
\begin{eqnarray}
\hat{X}_n(\tb{x}_0)
&=&
\hat{x}^{C^{\scriptscriptstyle{(k)}}}_{n-1}(\tb{x}_0)+\Delta\hat{W}_{n-1}^{\,\scriptscriptstyle{(k)}_{\vphantom{a_1}}}(\tb{x}_0) \nn\\
&&+h\sum_{|\boldsymbol{r}|=0}^N\frac{1}{\boldsymbol{r}!}\frac{\partial^{|\boldsymbol{r}|}}{\partial \tb{x}_{n-1}^{\boldsymbol{r}}}u^{\scriptscriptstyle{(k)}_{\vphantom{a_1}}}\big(\hat{\tb{x}}_{n-1}^C(\tb{x}_0),t_{n-1}\big)\Delta\hat{\tb{W}}_{n-1}^{\boldsymbol{r}}(\tb{x}_0)\nn\\
&&+\sum_{|\boldsymbol{r}|=0}^{N-1}\sum_{j=1}^d\frac{1}{\boldsymbol{r}!}\frac{\partial^{|\boldsymbol{r}|}}{\partial \tb{x}_{n-1}^{\boldsymbol{r}}}G^{\scriptscriptstyle{(k,j)}_{\vphantom{a_1}}}\big(\hat{\tb{x}}_{n-1}^C(\tb{x}_0),t_{n-1}\big)\Delta W^{{\,\scriptscriptstyle{(j)}_{\vphantom{a_1}}}}_{n}\Delta\hat{\tb{W}}_{n-1}^{\boldsymbol{r}}(\tb{x}_0) \nn\\
&&+R^{\scriptscriptstyle{(k)}}_{n,N}\big(\Delta\tb{W}_n,\Delta\hat{\tb{W}}_{n-1}(\tb{x}_0)\big), \label{eqJ3}
\end{eqnarray}
where the new remainder term is 
\begin{eqnarray}
R^{(k)}_{n,N}\big(\Delta\tb{W}_n,\Delta\hat{\tb{W}}_{n-1}(\tb{x}_0)\big) &\coloneq& hR^{u^{\scriptscriptstyle{(k)}}}_{n-1,N}\big(\Delta\hat{\tb{W}}_{n-1}(\tb{x}_0)\big) \nn\\
&& +\sum_{j=1}^d\Delta W^{{\,\scriptscriptstyle{(j)}_{\vphantom{a_1}}}}_{n}R^{G^{\scriptscriptstyle{(k,j)}}}_{n-1,N-1}\big(\Delta\hat{\tb{W}}_{n-1}(\tb{x}_0)\big). \nn
\end{eqnarray}
If we decompose $\hat{X}^{\vphantom{a}^{\scriptscriptstyle{(k)}}}_n(\tb{x}_0) = \hat{x}^{C^{\scriptscriptstyle{(k)}}}_{n}(\tb{x}_0)+\Delta\hat{W}_{n}^{\,\scriptscriptstyle{(k)}_{\vphantom{a_1}}}(\tb{x}_0)$ and then substitute
$$
\hat{\tb{x}}_n^C(\tb{x}_0)=\hat{\tb{x}}_{n-1}^C(\tb{x}_0)+h\boldsymbol{u}\big(\boldsymbol{\hat{x}}_{n-1}^C(\tb{x}_0,t_{n-1}\big),
$$
on the left-hand side of Eq. \eqref{eqJ3}, then we arrive at the identity \eqref{eqJ-1} in the statement of Theorem \ref{eff_noise_approx}. The convergence condition of the expansion is straightforward from Definition \ref{Conv_reg_def}.
\vspace{-0.0cm}\hfill$\blacksquare$
}

\cblack{
Let us remark that the polynomial approximation given in Theorem \ref{eff_noise_approx} can be written as a polynomial \textit{exclusively dependent} on the subsequence of independent noise increments $\Delta\tb{W}_m$, for $m=1,\ldots,n$, i.e., it is possible to write
\begin{equation}\label{W_span_only_noises}
\Delta\hat{W}_{n}^{\,\scriptscriptstyle{(k)}_{\vphantom{a_1}}}=\textrm{pol}_N\big(\Delta \tb{W}_n,\ldots,\Delta\tb{W}_{1}\big)+R^{\scriptscriptstyle{(k)}}_{n,N}\big(\Delta\tb{W}_n,\ldots,\Delta\tb{W}_{1}\big),
\end{equation}
where $\textrm{pol}_N(\cdots)$ denotes polynomial of order $N$. This fact can be easily verified by induction. specifically, expression \eqref{W_span_only_noises} shows that the convergence of the expansion at step $n$ depends only on the noise increments $\Delta \tb{W}_m$, $m=1,\ldots,n$. Therefore, in order to apply Theorem \ref{eff_noise_approx} in the analysis of Algorithm \ref{alg_main_ps_1}, we need to establish the conditions that $\Delta \tb{W}_m$ should satisfy in order to guarantee the convergence of the polynomial expansions of $\hat{\tb{X}}_n$ or $\Delta \hat{\tb{W}}_n(\tb{x}_0)$. 
}

\cblack{
From Eq.\eqref{W_span_only_noises} it can be seen that if the increments of the original noise process, $\Delta\tb{W}_{m}$, $m=1,\ldots,n$, are bounded, then the increments of the effective noise process, $\Delta\hat{\tb{W}}_{n}$, are bounded too. Lemma \ref{Bound_noise} below yields explicit bounds for the effective noise $\Delta\hat{\tb{W}}_{n}$ in terms of any available bound on $\Delta\tb{W}_{m}$.
}

\cblack{
\begin{lemma}\label{Bound_noise}
If there are finite constants $A^{\scriptscriptstyle{(j)}_{\vphantom{a_1}}}_n$ such that $\big|\Delta W^{\scriptscriptstyle{(j)}_{\vphantom{a_1}}}_n\big| < A^{\scriptscriptstyle{(j)}_{\vphantom{a_1}}}_n$ for every $n \ge 1$ and $j=1,\ldots,d$, then the constants recursively computed as
\begin{eqnarray}
\widetilde{A}^{\scriptscriptstyle{(k)}_{\vphantom{a_1}}}_{n,N}(\tb{x}_0) &=& \widetilde{A}^{\scriptscriptstyle{(k)}_{\vphantom{a_1}}}_{n-1,N}(\tb{x}_0)+h\!\!\sum_{|\boldsymbol{r}|=1}^N\frac{1}{\boldsymbol{r}!}\left|\frac{\partial^{|\boldsymbol{r}|}}{\partial \tb{x}_{n-1}^{\boldsymbol{r}}}u^{\scriptscriptstyle{(k)}_{\vphantom{a_1}}}\big(\hat{\tb{x}}_{n-1}^C(\tb{x}_0),t_{n-1}\big)\right|\widetilde{\tb{A}}_{n-1,N}^{\boldsymbol{r}}(\tb{x}_0) \nn \\
&& +\sum_{|\boldsymbol{r}|=0}^{N-1}\sum_{j=1}^d\frac{1}{\boldsymbol{r}!}\left|\frac{\partial^{|\boldsymbol{r}|}}{\partial \tb{x}_{n-1}^{\boldsymbol{r}}}G^{\scriptscriptstyle{(kj)}_{\vphantom{a_1}}}\big(\hat{\tb{x}}_{n-1}^C(\tb{x}_0),t_{n-1}\big)\right|A^{\scriptscriptstyle{(j)}_{\vphantom{a_1}}}_n\widetilde{\tb{A}}_{n-1,N}^{\boldsymbol{r}}(\tb{x}_0), \nn
\end{eqnarray}
where $\widetilde{\tb{A}}_{0,N}(\tb{x}_0)=0$, are finite and satisfy the inequalities
$$
\big|\Delta\hat{W}^{(k)}_{n,N} \big| < \widetilde{A}^{(k)}_{n,N}
$$ 
for every $n \ge 1$ and $k=1,\ldots,d$. Moreover, if 
$$
\widetilde{\tb{A}}_{n-1,N}(\tb{x}_0) \in \rho_{u^{(k)}}\big(\hat{\tb{x}}^C_{n-1}(\tb{x}_0),t_{n-1}\big) \cap_{j=1}^d \rho_{G^{(k,j)}}\big(\hat{\tb{x}}^C_{n-1}(\tb{x}_0),t_{n-1}\big),
$$ 
then
$$
\lim_{N\rightarrow\infty}R^{\scriptscriptstyle{(k)}}_{n,N}\big(\Delta\tb{W}_n,\Delta\hat{\tb{W}}_{n-1}(\tb{x}_0)\big)=0.
$$  
\end{lemma}
{\sc Proof:} It is straightforward from Theorem \ref{eff_noise_approx}.
}

\hfill$\blacksquare$

\cblack{
We can now apply the results above to provide a convergence condition for the recursive approximation of moments in Eq. \eqref{cond_moms_eff_noise}.
}

\cblack{
\begin{theorem}\label{expectation_sol}
Assume that the functions $\boldsymbol{u}$ and $\boldsymbol{G}$ in Eq. \eqref{dyn_sde} are real and smooth. For any positive integers $N$ and $n$, and any fixed initial condition $\tb{X}_0=\tb{x}_0$, we have the identity
\begin{eqnarray}
\mathds{E}\big[\Delta\hat{\tb{W}}^{\boldsymbol{r}}_{n}(\tb{x}_0)\big] &=& \sum_{|\boldsymbol{s}|+|\boldsymbol{r'}|=1}^N \!\!\!a_{\boldsymbol{r},n,N}^{(\boldsymbol{s},\boldsymbol{r}')}\big(\hat{\tb{x}}_{n-1}^C(\tb{x}_0),t_{n-1}\big)\mathds{E}\big[\Delta\tb{W}_{n}^{\,\boldsymbol{s}_{\vphantom{a_1}}}\big]\mathds{E}\big[\Delta\hat{\tb{W}}_{n-1,N}^{\boldsymbol{r'}}(\tb{x}_0)\big] \nn \\
&&+\mathds{E}\big[\boldsymbol{R}^{(\boldsymbol{r})}_{n,N}\big(\Delta\tb{W}_n,\Delta\hat{\tb{W}}_{n-1}(\tb{x}_0)\big)\big],
\label{eqJ4}
\end{eqnarray}
where $\boldsymbol{R}^{(\boldsymbol{r})}_{n,N}$ is the remainder term of the expansion and $a_{\boldsymbol{r},n,N}^{(\boldsymbol{s},\boldsymbol{r'})}\big(\hat{\tb{x}}_{n-1}^C(\tb{x}_0),t_{n-1}\big)$ are the constant coefficients of the expansion of the effective noise in Theorem \ref{eff_noise_approx}. Moreover, if there are finite constants $A^{\scriptscriptstyle{(j)}_{\vphantom{a_1}}}_n$ such that $\big|\Delta W^{\scriptscriptstyle{(j)}_{\vphantom{a_1}}}_n\big| < A^{\scriptscriptstyle{(j)}_{\vphantom{a_1}}}_n$ for every $n \ge 1$ and $j=1,\ldots,d$, then
\begin{equation}
\lim_{N\rightarrow\infty}\mathds{E}\big[\boldsymbol{R}^{(\boldsymbol{r})}_{n,N}\big(\Delta\tb{W}_n,\Delta\hat{\tb{W}}_{n-1}(\tb{x}_0)\big)\big]=0.
\label{eqJ4.5}
\end{equation}
\end{theorem}
}

\cblack{
{\sc Proof:} Note that the effective noise monomial $\Delta\hat{\tb{W}}^{\boldsymbol{r}}_{n}(\tb{x}_0)$ can be written as 
$$
\Delta\hat{\tb{W}}^{\boldsymbol{r}}_{n}(\tb{x}_0) = 
\prod_{i} \Delta\hat{W}_n^{(r_i)}(\tb{x}_0),
$$
where the factors $\Delta\hat{W}_n^{(r_i)}(\tb{x}_0)$ are expanded using Theorem \ref{eff_noise_approx} and then truncated to order $N$. We arrive at the identity \eqref{eqJ4}, after straightforward manipulations, by taking expectations and realising that 
$$
\mathds{E}\big[\Delta\tb{W}_{n}^{\,\boldsymbol{s}_{\vphantom{a_1}}}\Delta\hat{\tb{W}}_{n-1}^{\boldsymbol{r'}}(\tb{x}_0)\big]=\mathds{E}\big[\Delta\tb{W}_{n}^{\,\boldsymbol{s}_{\vphantom{a_1}}}\big]\mathds{E}\big[\Delta\hat{\tb{W}}_{n-1}^{\boldsymbol{r'}}(\tb{x}_0)\big],
$$ 
which is a consequence of Eq. \eqref{W_span_only_noises} and the independence of the noise increments.
\\
As for the convergence of the expansion \eqref{eqJ4}, Lemma \ref{Bound_noise} yields
$$
\lim_{N\rightarrow\infty}R^{\scriptscriptstyle{(k)}}_{n,N}\big(\Delta\tb{W}_n,\Delta\hat{\tb{W}}_{n-1}(\tb{x}_0)\big)=0,\qquad k=1,\ldots, v
$$
and, since the remainder term of the expansion of $\Delta\hat{\tb{W}}^{\boldsymbol{r}}_{n}(\tb{x}_0)$ is the addition of a finite sum of products involving 
$R^{\scriptscriptstyle{(k)}}_{n,N}\big(\Delta\tb{W}_n,\Delta\hat{\tb{W}}_{n-1}(\tb{x}_0)\big)$ for all $k=1,\ldots,v$, we obtain that
\begin{equation}
\lim_{N\rightarrow\infty}\boldsymbol{R}^{(\boldsymbol{r})}_{n,N}\big(\Delta\tb{W}_n,\Delta\hat{\tb{W}}_{n-1}(\tb{x}_0)\big)=0.
\label{eqJ5}
\end{equation}
Finally, if we take the expectation of \eqref{eqJ5} and apply the dominated convergence Theorem \cite{Re80} we arrive at \eqref{eqJ4.5} and complete the proof.
\\
\hfill$\blacksquare$
}

\cblack{
The results we have obtained are useful to guarantee convergence when the support of the dynamical noise $\tb{W}_n$ is bounded but, in general, this is not the case. However, even if in the most common models (Gaussian distributions, Gamma distributions, etc.) the support is not actually bounded, when the tails of a distribution decrease rapidly enough the support can be treated as bounded for numerical purposes. For example, if $W_n^{(k)} \sim \mathcal{N}(0,\sigma)$ then $\mathds{P}(|W_n^{(k)}|<3\sigma) > 0.9973$, i.e., $W_n^{(k)}$ is bounded with high probability. 
}

\cblack{
For a prescribed probability $P \in (0,1)$, let us choose the quantities
\begin{equation}\label{A_par_conv}
A^{\scriptscriptstyle{(j)}_{\vphantom{a_1}}}_n(P)\coloneq \inf \Big\{a\in\mathds{R}_0^{+}:\mathds{P}\big(\big|W_n^{\,\scriptscriptstyle{(j)}_{\vphantom{a_1}}}\big| < a\big) > P\Big\},\qquad j=1,\ldots,d,
\end{equation}
i.e., $A^{\scriptscriptstyle{(j)}_{\vphantom{a_1}}}_n(P)$ is an upper bound for $\big|W_n^{\,\scriptscriptstyle{(j)}_{\vphantom{a_1}}}\big|$ with probability $P$. One can combine bounds that hold with some probability $P$ and Lemma \ref{Bound_noise} to assess the convergence of the polynomial expansions of the moments of the effective noise in Algorithm \ref{alg_main_ps_1}.
}



\section{Approximation of 1-dimensional marginal densities}\label{pdf_sol_noise}

\cblack{
In this section we prove that the approximate 1-dimensional pdf's computed using Algorithm \ref{alg_main_ps_2} converge as the order of the Gram-Charlier expansion, $N$, increases, provided that the initial condition is fixed, $\tb{X}_0=\tb{x}_0$. When the initial condition is random, we further extend the latter result with the convergence of the Monte Carlo estimator in Eq. \eqref{law_tot_prob} as the number of samples $N_s^\prime$ increases.
}

\cblack{
\begin{theorem}\label{theorem_Charlier}
Let $X$ be a real random variable with pdf $f_{X}$ and characteristic function $\Psi_X$; then choose an auxiliary random variable $Z_\varphi$ with smooth pdf $\varphi$ and characteristic function $\Psi_\varphi$ such that $|\frac{\Psi_{X}(t)}{\Psi_{\varphi}(t)}|<\infty$ for all finite $t$. The density $f_{X}$ can be expanded with respect to the derivatives of $\varphi$ as
\beq
f_{X}(x)=\left[1+\sum_{r=1}^{N}\frac{(-1)^{r}}{r!}C_{r}\big[X,Z_\varphi\big]\frac{\diff^{r}}{\diff {x}^{r}}\right]\!\varphi(x)+R_{N}\big(f_{X},\varphi;x\,\big),
\label{eqX00}
\eeq
where $C_{r}\big[X,Z_\varphi\big]$ are the coefficients of the expansion defined in Eq.\,\eqref{Coeff_GC} and $R_{N}\big(f_{X},$ $\varphi;x\,\big)$ is a remainder term.
\end{theorem}
}


\cblack{
{\sc Proof:} 
We write $\Psi_{X}$ as
\beqa
\Psi_{X}(t) &=& \frac{\Psi_{X}(t)}{\Psi_{\varphi}(t)}\Psi_{\varphi}(t) \nn \\
&=& \exp\Big(\log\big(\Psi_{X}(t)\big)-\log\big(\Psi_{\varphi}(t)\big)\Big)\Psi_{\varphi}(t) \nn \\
&=& \left[\exp\left(\sum_{r=0}^N\frac{\kappa_{r}(X)-\kappa_{r}(Z_\varphi)}{r!}(it)^{r}\right)+\textrm{R}_N\left(\frac{\Psi_{X}}{\Psi_{\varphi}};t\right)\right]\!\Psi_{\varphi}(t),
\label{eqX01}
\eeqa
where $\textrm{R}_N\left(\frac{\Psi_{X}}{\Psi_{\varphi}};t\right)$ is the remainder of the Taylor expansion of function $\Psi_{X}/\Psi_{\varphi}$. 
If we expand the exponential function in \eqref{eqX01} in terms of Bell polynomials \cite{Ab05} and then compute the inverse Fourier transform on both sides of the equation 
we arrive at
$$
f_{X}(x)=\left[1+\sum_{r=1}^{N}\frac{(-1)^{r}}{r!}C_{r}\big[X,Z_\varphi\big]\frac{\diff^{r}}{\diff {x}^{r}}\right]\!\varphi(x) +\frac{1}{2\pi}\!\!\int\limits_{\;\;\;\,\mathds{R}}\textrm{R}_N\bigg(\frac{\Psi_{X}}{\Psi_{\varphi}};t\bigg)\Psi_{\varphi}(t)\e^{-ixt}\diff t,
$$
where the second term on the r.h.s. is the remainder in Eq. \eqref{eqX00}. \hfill$\blacksquare$
}

\cblack{
Many families of orthogonal polynomials are related to specific probability distributions \cite{Ch11} in the sense that there are formulas to generate the polynomials from the derivatives of probability densities (the so-called Rodrigues formulas \cite{Ra81}). In particular, the class of probabilistic Hermite polynomials are orthogonal w.r.t. the Gaussian distribution and the Rodrigues formula for them is given by Eq. \eqref{Rodrigues}.
}

\cblack{
If we let the auxiliary pdf $\varphi$ be Gaussian distribution, the Gram--Charlier expansion of $f_{X}$ in Theorem \ref{theorem_Charlier} reduces to a series of Hermite polynomials multiplied by $\varphi$. Hence, the convergence of expression \eqref{eqX00} becomes a standard problem, similar to the convergence of the PCE \eqref{pol_exp_1} in Section \ref{PCE_subsec}. Indeed, if $\varphi$ is Gaussian, the series in \eqref{eqX00} is termed a Gram--Charlier expansion of type A and it can be expected to converge when\footnote{We construct the class of real $L^2$ functions w.r.t. a density $\varphi : S \mapsto (0,\infty)$ as 
$$
L^2(S,\varphi) := \left\{
	h : S \mapsto \mathds{R} \text{ such that } \int_S h(x)\varphi(x)\diff x < \infty
\right\}. 
$$
}
$
\frac{f_{X}}{\varphi}\in L^2\left(\mathds{R},\varphi\right)
$ (see \cite{Ch11}). 
}

\cblack{
In the sequel, we restrict our attention to the Gram-Charlier expansion of type A and hence assume that the auxiliary pdf $\varphi$ used to approximate the $k$-th 1-dimensional marginal pdf $f_{\hat{X}_n^{(k)} | \tb{x}_0}$ is Gaussian, with mean $\mu_n^{(k)}(\tb{x}_0)$ and standard deviation $\sigma_n^{(k)}(\tb{x}_0)$. We specifically denote it as $\varphi_{n,\tb{x}_0}^{(k)}(x)$ (note the dependence on the initial condition $\tb{x}_0$).
}

\cblack{
Below, we establish some regularity assumptions and then use them to provide an explicit convergence theorem for the approximations of $f_{\hat{X}_n^{(k)} | \tb{x}_0}$.
}


\cblack{
\begin{assumption}\label{ass-f-bounds}
Let $\mathrm{supp}(f)$ denote the support of function $f$ and let $f_{\Delta\tb{W}_m}$ denote the pdf of the random vector of noise increments at time $m$, $\tb{W}_m$. There are bounded sets $D_m \subset \mathds{R}^d$, $m=1, \ldots, n$, such that 
$$
\mathrm{supp}(f_{\Delta\tb{W}_m}) \subseteq D_m.
$$
Moreover, there is a sequence of finite constants $M_m$, $m = 1, \ldots, n$, that satisfy the inequalities
$$
\sup_{x\,\in\, \mathds{R} }\big(f_{\Delta\tb{W}_m}(x)\big)\leq M_m.
$$
\end{assumption}
}

\cblack{
\begin{assumption} \label{ass-a-bounds}
There are finite constants $\big\{A^{\scriptscriptstyle{(j)}_{\vphantom{a_1}}}_m : m=1, \ldots, n; \quad j=1, \ldots, v \big\}$ such that $| \Delta W_m^{(j)} | < A^{(j)}_m$, for $1 \le j \le v$ and $1\le m \le n$, and 
\beq
\widetilde{\tb{A}}_{m,N}(\tb{x}_0) \in 
\rho_{u^{(k)}}\big(\hat{\tb{x}}^C_{m}(\tb{x}_0),t_{m}\big)
\cap_{j=1}^d
\rho_{G^{(k,j)}}\big(\hat{\tb{x}}^C_{m}(\tb{x}_0),t_{m}\big)
\eeq
where $k$-th entry of the $v$-dimensional vector $\widetilde{\tb{A}}_{m,N}(\tb{x}_0)$ is constructed as
\beqa
 \widetilde{A}^{\scriptscriptstyle{(k)}_{\vphantom{a_1}}}_{m,N}(\tb{x}_0) &=& \widetilde{A}^{\scriptscriptstyle{(k)}_{\vphantom{a_1}}}_{m-1,N}(\tb{x}_0)+h\!\!\sum_{|\boldsymbol{r}|=1}^N\frac{1}{\boldsymbol{r}!}\left|\frac{\partial^{|\boldsymbol{r}|}}{\partial \tb{x}_{m-1}^{\boldsymbol{r}}}u^{\scriptscriptstyle{(k)}_{\vphantom{a_1}}}\big(\hat{\tb{x}}_{m-1}^C(\tb{x}_0),t_{m-1}\big)\right|
\widetilde{\tb{A}}_{m-1,N}^{\boldsymbol{r}}(\tb{x}_0) \nn \\
&& +\sum_{|\boldsymbol{r}|=0}^{N-1}\sum_{j=1}^d\frac{1}{\boldsymbol{r}!}\left|\frac{\partial^{|\boldsymbol{r}|}}{\partial \tb{x}_{m-1}^{\boldsymbol{r}}}G^{\scriptscriptstyle{(k,j)}}\big(\hat{\tb{x}}_{m-1}^C(\tb{x}_0),t_{m-1}\big)\right|
A^{\scriptscriptstyle{(j)}_{\vphantom{a_1}}}_m\widetilde{\tb{A}}_{m-1,N}^{\boldsymbol{r}}(\tb{x}_0),
\eeqa
with initial condition $\widetilde{\tb{A}}_{0,N}(\tb{x}_0)=0$.
\end{assumption}
}


\cblack{
Assumption \ref{ass-f-bounds} states that the support of the noise components is bounded, while Assumption \ref{ass-a-bounds} guarantees that finite noise increments $\Delta W_m^{(k)}$ yield finite effective noise terms $\Delta \hat{W}^{(k)}_m$ and enables us to apply Lemma \ref{Bound_noise}. Given the above regularity assumptions we can provide guarantees on the approximation of the marginal densities $f_{\hat{X}_n^{(k)}|\tb{x}_0}(x)$.
}

\cblack{
\begin{theorem}\label{Thm_Charlier_conv}
Let the functions $\boldsymbol{u}$ and $\boldsymbol{G}$ in the SDE \eqref{dyn_sde} be smooth, let Assumptions \ref{ass-f-bounds} and \ref{ass-a-bounds} hold and let $\tb{x}_0$ be a fixed initial condition. Then, the type A Gram-Charlier expansion of the 1-dimensional marginal pdf of $\hat{X}_n^{(k)}(\tb{x}_0)$, $k \in \{1, \ldots, v \}$, can be written as
\beqa
f_{\hat{X}^{\scriptscriptstyle{(k)}}_{n}|\tb{x}_0}(x)&=&\left[1+\sum_{r=1}^N\frac{1}{r!\,\sigma^{\scriptscriptstyle{(k)}}_{n}(\tb{x}_0)^{\vphantom{A}^{r}}}C_{r}\Big[\hat{X}^{(k)}_{n,N}(\tb{x}_0),Z^{(k)}\!(\tb{x}_0)\Big]\vphantom{H_{e_r}\left(\frac{x-\mu^{(k)}_{n}(\tb{x}_0)}{\sigma^{(k)}_{n}(\tb{x}_0)}\right)} H_{e_r}\left(\frac{x-\mu^{(k)}_{n}(\tb{x}_0)}{\sigma^{(k)}_{n}(\tb{x}_0)}\right)\right]\times \nn \\
&&\times \varphi^{(k)}_{n,\tb{x}_0}\!(x)+R^{(k)}_{n,N\vphantom{N}}\big(f_{\hat{X}^{(k)}_{n}|\tb{x}_0},\varphi^{(k)}_{n,\tb{x}_0}(\cdot);x\,\big), \nn
\eeqa
where $\big\{H_{e_r}\big\}_{{r=0}}^\infty$ are the probabilistic Hermite polynomials given by Eq.\,\eqref{Rodrigues}, $Z^{(k)}(\tb{x}_0)$ is a random variable with pdf $\varphi_{n,\tb{x}_0}^{(k)}$ and the remainder term vanishes as the truncation order $N$ is increased, i.e.,
\beq
\lim_{N\rightarrow\infty} R^{(k)}_{n,N\vphantom{N}}\big(f_{\hat{X}^{(k)}_{n}|\tb{x}_0},\varphi^{(k)}_{n,\tb{x}_0}(\cdot);x\,\big) =0.
\label{eqXX10}
\eeq
\end{theorem}
}

{\sc Proof:} 
\cblack{
The type A Gram--Charlier expansion of $f_{\hat{X}^{(k)}_{n}|\tb{x}_0}(x)$ is immediately obtained from Eq. \eqref{char_gen_form} when the auxiliary pdf is Gaussian (namely, $\varphi = \varphi^{(k)}_{n,\tb{x}_0}$). Additionally, we need to prove that Eq. \eqref{eqXX10} holds, which takes more effort. Specifically, hereafter we prove that the function $\frac{ f_{\hat{X}^{(k)}_{n}|\tb{x}_0} }{ \varphi^{(k)}_{n,\tb{x}_0} }$ belongs to $L^2(\mathds{R},\varphi^{(k)}_{n,\tb{x}_0})$, which, in turn, implies that $R^{(k)}_{n,N\vphantom{N}}\big(f_{\hat{X}^{(k)}_{n}|\tb{x}_0},\varphi^{(k)}_{n,\tb{x}_0}(\cdot);x\,\big) \stackrel{N\rightarrow\infty}{\longrightarrow} 0$ in $L^2$ (see \cite{Ch11}).
}

\cblack{
First, we prove using an induction argument that the pdf $f_{\hat{\tb{X}}_{n-1}|\tb{x}_0}$ is bounded and it has a bounded support. Let us assume that at time $n-1$ there are a bounded set $\hat{D}_{n-1} \subset \mathds{R}^v$ and a finite constant $\hat{M}_{n-1}$ such that
\beq
\mathrm{supp}\big(f_{\hat{\tb{X}}_{n-1}|\tb{x}_0}\big) \subseteq \hat{D}_{n-1}\!\subset\mathds{R}^{v}\qquad\mbox{and}\qquad\sup_{x\,\in\,\mathds{R}}\big(f_{\hat{\tb{X}}_{n-1}|\tb{x}_0}(x) \big)\leq\hat{M}_{n-1} < \infty. 
\label{eqXX20}
\eeq
From the expression of the Euler--Maruyama integrator in Eq. \eqref{sto_integrator}, we can write the pdf of $\hat{\tb{X}}_n$ in terms of the densities of $\hat{\tb{X}}_{n-1}$ and $\Delta\tb{W}_n$ as
\beqa
f_{\hat{\tb{X}}_{n}|\tb{x}_0}(\tb{x}_n) &=& 
\int_{\mathds{R}^{v+d}} f_{\hat{\tb{X}}_{n-1}|\tb{x}_0}(\tb{x}_{n-1}) f_{\Delta{\tb{W}}_{n}}(\tb{w}_{n}) \times \nn \\
&& \times \delta\left(
	\tb{x}_n-\tb{x}_{n-1} - h\tb{u}\left(
		\tb{x}_{n-1},t_{n-1}
	\right) - \tb{G}\left(
		\tb{x}_{n-1},t_{n-1}
	\right) \tb{w}_{n}
\right) \diff\tb{x}_{n-1}\diff\tb{w}_{n},
\label{eqXX21}
\eeqa
where $\delta(\cdot)$ denotes the Dirac delta function (see Eq. {\sc 4.34} in \cite{Re98}). Using Assumption \ref{ass-f-bounds} and the induction hypothesis \eqref{eqXX20} we obtain an upper bound for the pdf $f_{\hat{\tb{X}}_{n}|\tb{x}_0}(\tb{x}_n)$ in Eq. \eqref{eqXX21} of the form
\beqa
f_{\hat{\tb{X}}_{n}|\tb{x}_0}(\tb{x}_n) &\leq& \nn \\
\hat{M}_{n-1}M_n \int_{\hat{D}_{n-1}\times D_n} \delta\Big(\tb{x}_n-\tb{x}_{n-1}-h\tb{u}\big(\tb{x}_{n-1},t_{n-1}\big) - \tb{G}\big(\tb{x}_{n-1},t_{n-1}\big)\tb{w}_{n}\Big)\diff\tb{x}_{n-1}\diff\tb{w}_{n} &\le& \nn\\
\hat{M}_{n-1}M_n, && \nn
\eeqa
hence 
\beq
\sup_{x \in \mathds{R}}\left( 
	f_{\hat{\tb{X}}_n | \tb{x}_0}(x) 
\right) \leq \hat{M}_n < \infty, 
\quad \text{where} \quad
\hat{M}_n = \hat{M}_{n-1} M_n.
\label{eqXX22}
\eeq
Moreover, since $\tb{u}$ and $\tb{G}$ are smooth and $\hat{D}_{n-1} \times D_n$ is bounded, all the solutions of the equation
$$
\tb{x}_n-\tb{x}_{n-1}-h\tb{u}\big(\tb{x}_{n-1},t_{n-1}\big)-\tb{G}\big(\tb{x}_{n-1},t_{n-1}\big)\tb{w}_{n}=0
$$
necessarily lie in a bounded set $\hat{D}_n \subset \mathds{R}^v$, which implies that
\beq
\mathrm{supp}\big(f_{\hat{\tb{X}}_n|\tb{x}_0}\big) \subseteq \hat{D}_n \subset \mathds{R}^{v}.
\label{eqXX23}
\eeq
To complete the induction argument, we need to prove that
$$
\mathrm{supp}\big(f_{\hat{\tb{X}}_{1}|\tb{x}_0}\big) \subseteq \hat{D}_{1}\!\subset\mathds{R}^{v}\qquad\mbox{and}\qquad\sup_{x\,\in\,\mathds{R}}\big(f_{\hat{\tb{X}}_{1}|\tb{x}_0}(x)\big)\leq\hat{M}_{1}<\infty.
$$
for some bounded set $\hat{D}_1$ and some finite constant $\hat{M}_1$. Resorting again to the expression of the Euler--Maruyama integrator \eqref{sto_integrator} and Assumption \ref{ass-f-bounds} we obtain the inequalities
\beqa
f_{\hat{\tb{X}}_{1}|\tb{x}_0}(\tb{x}_1) &=& 
\int_{\mathds{R}^{d}} f_{\Delta{\tb{W}}_{1}}(\tb{w}_{1})\delta\Big(\tb{x}_1-\tb{x}_{0} - h\tb{u}\big(\tb{x}_{0},t_{0}\big) - \tb{G}\big(\tb{x}_{0},t_{0}\big)\tb{w}_{1}\Big)\diff\tb{w}_{1} \nn \\
&\leq& M_1\!\!\!\!\int\limits_{\quad\,\mathds{R}^{d}}\!\!\delta\Big(\tb{x}_1-\tb{x}_{0} - h\tb{u}\big(\tb{x}_{0},t_{0}\big)-\tb{G}\big(\tb{x}_{0},t_{0}\big)\tb{w}_{1}\Big)\diff\tb{w}_{1} \nn \\
&\le& M_1 < \infty, \nn
\eeqa
hence $\hat{M}_1 = M_1$ and, by the same reasoning as in the induction step, the solutions of the equation $\tb{x}_1-\tb{x}_{0} - h\tb{u}\big(\tb{x}_{0},t_{0}\big)-\tb{G}\big(\tb{x}_{0},t_{0}\big)\tb{w}_{1} = 0$ lie in a bounded set $\hat{D}_1 \subset \mathds{R}^v$ which contains the support of $f_{\hat{\tb{X}}_{1}|\tb{x}_0}$.
}

\cblack{
The bounds in \eqref{eqXX22} and \eqref{eqXX23} imply that $f_{\hat{\tb{X}}_n|\tb{x}_0} \in L^2\left( \mathds{R},\varphi_{n,\tb{x}_0}^{(k)}(x) \diff x \right)$. To show it, let $\hat{D}_n^{(k)}$ be the projection of the bounded set $\hat{D}_n \subset \mathds{R}^v$ along the $k$-th dimension and then note that 
\beq
\int_{\mathds{R}} 
\left(
	\frac{
		f_{\hat{X}^{(k)}_{n}|\tb{x}_0}(x)
	}{
		\varphi^{(k)}_{n,\tb{x}_0}(x)
	}
a
\right)^2 \varphi^{(k)}_{n,\tb{x}_0}(x) \diff x \leq
\frac{
	\int_{\hat{D}^{(k)}_n} \left(
		f_{\hat{X}^{(k)}_{n}|\tb{x}_0}(x)
	\right)^2 \diff x
}{
	\inf_{x \in \hat{D}^{(k)}_n} \varphi^{(k)}_{n,\tb{x}_0}(x)
}<\infty,
\label{eqXX24}
\eeq
where the second inequality holds because 
\begin{itemize}
\item $\inf_{x \in \hat{D}^{(k)}_n} \varphi^{(k)}_{n,\tb{x}_0}(x) > 0$, since $\varphi^{(k)}_{n,\tb{x}_0}$ is Gaussian and $\hat{D}^{(k)}_n$ is bounded, and
\item the marginal density $f_{\hat{X}^{(k)}_{n}|\tb{x}_0}$ is bounded because the joint density $f_{\hat{\tb{X}}_{n}|\tb{x}_0}$ is bounded.
\end{itemize}
The inequality \eqref{eqXX24} yields $f_{\hat{\tb{X}}_n|\tb{x}_0} \in L^2\left( \mathds{R},\varphi_{n,\tb{x}_0}^{(k)}(x) \diff x \right)$ which, in turn, implies that \eqref{eqXX10} holds \cite{Ch11}.
\hfill$\blacksquare$
}

\cblack{
Theorem \ref{Thm_Charlier_conv} states that the estimates of the 1-dimensional marginal pdf's $f_{\hat{X}^{(k)}_n | \tb{x}_0}(x)$ converge pointwise, for any fixed $\tb{x}_0$ and $x \in D_n$, as the truncation order $N$ increases. When the initial condition is random, the natural estimate to compute is the Monte Carlo approximation in Eq. \eqref{law_tot_prob}. The proposition below guarantees that, under similar assumptions as in Theorem \ref{Thm_Charlier_conv}, the Monte Carlo estimator converges to 
$$
f_{\hat{X}^{(k)}_n}(x) = \mathds{E}\left[
	f_{\hat{X}^{(k)}_n | \tb{X}_0}(x)
\right]
$$ almost surely (a.s.) for any $x \in D_n$.
}

\cblack{
\begin{proposition}\label{MC_lemma}
Let $\tb{X}_{0,j}$, $j=1, \ldots, N_s'$, be i.i.d. samples form the initial distribution with pdf $f_{\tb{X}_0}$. If Assumptions \ref{ass-f-bounds} and \ref{ass-a-bounds} hold and $f_{\tb{X}_0}$ is bounded with bounded support, then
\beq
\lim_{N_s'\rightarrow\infty} \left[
	\lim_{N\rightarrow\infty} \frac{1}{N'_s}\sum_{j=1}^{N'_s}f_{\hat{X}^{(k)}_{n,N}|\tb{X}'_{0,j}}(x)  
\right] = f_{\hat{X}^{(k)}_{n}}(x) \quad \text{a.s.},
\label{eqSLLN}
\eeq
for $k=1,\ldots,v$.
\end{proposition}
}

\cblack{
{\sc Proof:} For any $N_s'\in\mathds{N}$, Theorem \ref{Thm_Charlier_conv} yields
$$
\lim_{N\rightarrow\infty} \frac{1}{N'_s}\sum_{j=1}^{N'_s}f_{\hat{X}^{(k)}_{n,N}|\tb{X}'_{0,j}}(x) = \frac{1}{N'_s}\sum_{j=1}^{N'_s}f_{\hat{X}^{(k)}_{n}|\tb{X}'_{0,j}}(x). 
$$
Moreover, since $f_{\tb{X}_0}$ is bounded and has a bounded support, the same argument as in the proof of Theorem \ref{Thm_Charlier_conv} shows that the pdf's $f_{\hat{X}^{(k)}_{n}|\tb{X}'_{0,j}}$ are uniformly bounded\footnote{The bounds $\hat{M}_n$ in the proof of Theorem \ref{Thm_Charlier_conv} depend on the initialization, i.e., $\hat{M}_n=\hat{M}_n(\tb{X}_0)$. However, the bounds $\hat{M}_n(\tb{X}_0)$ are continuous by construction and, since the support of $\tb{X}_0$ is bounded, $\sup_{\tb{X}_0} \hat{M}_n(\tb{X}_0) < \infty$.}, hence $\mathds{E}\left[ \left( f_{\hat{X}^{(k)}_{n}|\tb{X}'_{0,j}}(x) \right)^2 \right] < \infty$ and the strong law of large numbers yields Eq. \eqref{eqSLLN}.
\vspace{0.cm}\hfill$\blacksquare$
}


\section{Numerical examples}\label{Num_examples}

\cblack{
In order to illustrate the performance of the proposed uncertainty quantification scheme we provide two numerical examples. In both of them, we compare the solution of the dynamics of a Keplerian orbit in two-dimensional space, perturbed by an additive Wiener process, using Algorithm \ref{alg_main_ps_1} for the approximation of moments and a Monte Carlo simulation with $10^4$ samples as a baseline. The two examples differ essentially in the choice of initial condition, which is fixed for the first set of simulations while we assume it random (with a Gaussian distribution) for the second example.
}

\cblack{
We start from the general equation \eqref{dyn_sde}. The state vector of the orbiting object has dimension $d=4$ and we denote it as 
$
\tb{X}(t)=\left(
x(t),
y(t),
v_x(t),
v_y(t)
\right)^\top
$, where $\left(x(t),y(t)\right)$ is the object position in km and $\left(v_x(t),v_y(t)\right)$ is its velocity in km/s, respectively, in 2-dimensional space. For the Keplerian dynamics, the drift coefficient $\tb{u}(\tb{X},t)$ can be written as \cite{Va07}
\begin{equation}
\tb{u}(\tb{X},t)=\begin{pmatrix}
v_x(t)\\
v_y(t)\\
\displaystyle{\frac{-\mu x(t)}{\big[x(t)^2+y(t)^2\big]^{3/2}}}\\
\displaystyle{\frac{-\mu y(t)}{\big[x(t)^2+y(t)^2\big]^{3/2}}}\\
\end{pmatrix},
\end{equation}
where $\mu$ is the standard gravitational parameter, and we set the diffusion coefficient as the $4 \times 4$ diagonal matrix $\tb{G}(\tb{X},t)=\diag\left( 0, 0, \sigma_3, \sigma_4 \right)$, where $\sigma_3$ and $\sigma_4$ are known positive constants. The noise process $\tb{W}(t)$ is a standard $4 \times 1$ Wiener process. Physically, $\diff\tb{W}(t)$ represents a stochastic perturbation in the acceleration of the orbiting object. The numerical values used for the simulation are summarised in Table \ref{Table_par}.
}

\begin{table}[H]
\cblack{
\begin{center}
	\vspace{0.cm}\begin{tabular}{|m{2.2cm}|m{2cm}||m{6.7cm}|}
	\hline
	\hfil\textbf{Parameters} & \hfil\textbf{Value} & \hfil\textbf{Description}\\	
	\hline\hline
	\hfil$\mu$ & \hfil $3.986$ $\frac{\text{km}^3}{\text{s}^2}$ & Standard gravitational parameter.\\
	\hline
	\hfil$h$ & \hfil $0.1\,$s & Step-size for time discretisation.\\
	\hline
	\hfil$t_0$ & \hfil $0\,$days & Initial time.\\
	\hline
	\hfil$t_n$ & \hfil $1.5\,$days & Final time.\\
	\hline
	\hfil$n$ & \hfil$1,296,000$ & Number of discrete-time steps in the simulations, namely, $n=\lceil(t_n-t_0)/h\rceil$ where $\lceil\,\cdot\,\rceil$ denotes the ceiling function.\\
	\hline
	\hfil$N$ & \hfil $2$ & Order of the polynomial expansions.\\
	\hline
	\hfil$N_{\textrm{PCE}}$ & \hfil $6$ & Truncation order of PCE (for the second example only).\\
	\hline
	\hfil$N_s$ & \hfil $420$ & Number of samples (for the second example only).\\
	\hline
	\hfil$N'_s$ & \hfil $10^5$ & Number of samples of $\tb{X}_0$ generated to reconstruct the marginal pdf's $f_{\hat{X}^{(k)}_{n,N}}(x)$.\\
	\hline
	\hfil$\sigma_3, \sigma_4$ & \hfil $10^{-5}$ & Scale parameters in the diffusion term $\tb{G}(\tb{X},t)$.\\
	\hline
	\end{tabular}
\end{center}
\caption{Simulation parameters.}
\label{Table_par}
}
\end{table}

\cblack{
All computer experiments have been performd using Matlab R2018b running on a Mac Book Pro computer equipped with a $2.3$~GHz Intel Core i5 CPU and 16 GB of RAM.
}


\label{exper_1}\subsection{Fixed initial condition}

\cblack{
For the first experiment, we fix the initial condition as
\begin{equation}\label{in_cond_state}
\tb{x}_0=\begin{pmatrix}
200+R_T\\
0\\
0\\
\displaystyle{\sqrt{\frac{\mu}{200+R_T}}}
\end{pmatrix},
\end{equation}
where $R_T=6.378\times 10^3\,$km is the Earth radius. By taking a known initial condition, we can asses the moment and density approximations when the only source of uncertainty is the dynamical noise $\tb{W}(t)$ and, therefore, we avoid any PCE approximation.
}

\cblack{
The initial state $\tb{x}_0$ has been chosen to simulate the evolution of a circular orbit at 200 km above the Earth surface. At this low altitude, it is relevant to use a SDE to represent the orbital dynamics because the object motion depends on the atmosphere drag which, in turn, depends on several parameters (atmosphere density, mass, volume, shape of the object, etc.) which are often difficult to determine in practice \cite{Va07}. The diffusion term $\tb{G}(\tb{X},t) \diff \tb{W}(t)$ may account for these uncertainties.
}

\cblack{
Table \ref{Table_mean_fix} shows a comparison between the outcomes, at the final time $t_n$, of Algorithm \ref{alg_main_ps_1} and the baseline Monte Carlo method with $10^4$ independent trajectories generated using the Euler-Maruyama scheme \eqref{sto_integrator}. The first column in the table displays the expected values of $x$, $y$, $v_x$ and $v_y$ computed with Algorithm \ref{alg_main_ps_1}, while the second column shows the Monte Carlo estimates for each state variable. The third column displays the absolute differences between the first and second columns, and the fourth column shows the relative difference (with the Monte Carlo estimates taken as reference). We can observe that both methods yield very similar outputs, with the relative differences of just $\sim 0.1\%$ for $x$ and $v_y$, and $\sim 0.006\%$ for $y$ and $v_x$.  
}

\vspace{1cm}\begin{table}[htbp]
\cblack{
\begin{center}
\begin{tabular}{|c|c|c|c|c|}\hline
&Algorithm \ref{alg_main_ps_1}	&Monte Carlo, 	&Absolute  	&Relative\\
&								&$10^4$ samples	&difference	&difference\\
\hline\hline
$x$ & $\vphantom{K^{B^2}}$ $\phantom{-}1.18\times 10^3$\,km & $\phantom{-}1.18\times 10^3$\,km & $1.53$\,km & $1.30\times 10^{-3}$\\\hline
$y$ & $\vphantom{K^{B^2}}-6.66\times 10^3$\,km & $-6.66\times 10^3$\,km & $4.61\times 10^{-1}$\,km & $6.92\times 10^{-5}$\\\hline
$v_x$ & $\vphantom{K^{B^2}}\phantom{-}7.48$\,km/s & $7.48$\,km/s & $4.42\times 10^{-4}$\,km/s & $5.90\times 10^{-5}$\\\hline
$v_y$ & $\vphantom{K^{B^2}_{B_2}}\phantom{-}1.32$\,km/s & $1.33$\,km/s & $1.77\times 10^{-3}$\,km/s & $1.34\times 10^{-3}$\\\hline
\end{tabular}
\caption{\cblack{Estimate of $\mathds{E}[\tb{X}(t_n)]$ with the moment-computation Algorithm \ref{alg_main_ps_1}, compared with standard Monte Carlo estimates. The initial condition is fixed}}
\label{Table_mean_fix}
\end{center}
}
\end{table}

\cblack{
Table \ref{Table_Cov_fix} shows a comparison between the estimates of the second order moments of $\tb{X}(t_n)$ computed via Algorithm \ref{alg_main_ps_1} and the standard Monte Carlo method that runs the Euler-Maruyama scheme $10^4$ times. The first shows the covariance matrix of $\tb{X}(t_n)$ as output by Algorithm \ref{alg_main_ps_1}, while the second row shows the Monte Carlo estimate. The entry-wise absolute and relative differences between the two matrices are displayed in the third and fourth rows of the table, respectively. The differences are larger than for the first-order moments, yet the two methods still yield very similar outputs.  
}

\begin{table}[htb]
\cblack{
\begin{center}
\begin{tabular}{|c||c|}\hline
Algorithm \ref{alg_main_ps_1}&
$\vphantom{\begin{matrix}
\vspace{-0.3cm}a\\
a\\
a\\
a\\
a\\
a
\end{matrix}}\begin{pmatrix}
\phantom{-}5.92\times 10^{5} &  \phantom{-} 1.09\times 10^{5} &-1.16\times 10^{2} & \!\!\!\!\phantom{-} 6.67\times 10^{2}\\
\phantom{-}1.09\times 10^{5} & \phantom{-} 1.82\times 10^{4} &-1.91\times 10^{1} & \!\!\!\!\phantom{-}1.23\times 10^{2} \\
 -1.16\times 10^{2} & -1.91\times 10^{1} &\phantom{-}\;\;\,2.02\times 10^{-2} & -1.31\times 10^{-1}\\
\phantom{-}6.67\times 10^{2} & \phantom{-}1.23\times 10^{2} & \;\;\,-1.31\times 10^{-1} &\phantom{-}7.50\times 10^{-1}
\end{pmatrix}$\\
\hline
Monte Carlo
&
$\vphantom{\begin{matrix}
\vspace{-0.3cm}a\\
a\\
a\\
a\\
a\\
a
\end{matrix}}\begin{pmatrix}
 \phantom{-}5.78\times 10^{5} & \phantom{-}1.05\times 10^{5} & -1.11\times 10^{2} &  \!\!\!\!\phantom{-}6.51\times 10^{2}\\
 \phantom{-}1.05\times 10^{5} &  \phantom{-}2.30\times 10^{4} &  -2.46\times 10^{1} &  \!\!\!\!\phantom{-}1.18\times 10^{2}\\
 -1.11\times 10^{2} & -2.46\times 10^{1} & \;\;\,\phantom{-}2.64\times 10^{-2} &-1.25\times 10^{-1}\\
\phantom{-}6.51\times 10^{2} & \phantom{-}1.18\times 10^{2}  & \;\;\,-1.25\times 10^{-1}&\phantom{-}7.32\times 10^{-1}
\end{pmatrix}$\hspace{-0.14cm}\\\hline
$\begin{matrix}
\text{Absolute}\\
\text{differences}
\end{matrix}
$
& 
$\vphantom{\begin{matrix}
\vspace{-0.3cm}a\\
a\\
a\\
a\\
a\\
a
\end{matrix}}\hspace{-0.1cm}\begin{pmatrix}
 1.40\times 10^{4} &  4.56\times 10^{3} & 4.97 & \hspace{-0.2cm}1.59\times 10^{1}\\
 4.56\times 10^{3} &  4.82\times 10^{3} & 5.48 & 5.02 \\
 4.97 & 5.48 & 6.23\times 10^{-3} & 5.46\times 10^{-3}\\
 1.59\times 10^{1} & 5.02 & 5.46\times 10^{-3} & 1.80\times 10^{-2}
\end{pmatrix}$\hspace{-0.14cm}\\\hline
$\begin{matrix}
\text{Relative}\\
\text{differences}
\end{matrix}
$ & 
$\vphantom{\begin{matrix}
\vspace{-0.3cm}a\\
a\\
a\\
a\\
a\\
a
\end{matrix}}\hspace{-0.1cm}\begin{pmatrix}
2.43\times 10^{-2}& 4.35\times 10^{-2} & 4.47\times 10^{-2} & 2.44\times 10^{-2} \\
4.35\times 10^{-2}& 2.10\times 10^{-1} & 2.23\times 10^{-1} & 4.24\times 10^{-2} \\
4.47\times 10^{-2}& 2.23\times 10^{-1} & 2.36\times 10^{-1} & 4.36\times 10^{-2} \\
2.44\times 10^{-2}& 4.24\times 10^{-2} & 4.36\times 10^{-2} & 2.46\times 10^{-2} \\
\end{pmatrix}$\\\hline
\end{tabular}
\caption{\cblack{Estimates of the covariance matrix of $\tb{X}(t_n)$ computed via Algorithm \ref{alg_main_ps_1} and standard Monte Carlo, with $10^4$ independent samples, with fixed initial condition.}}
\label{Table_Cov_fix}
\end{center}
}
\end{table}

\cblack{
Since Algorithm \ref{alg_main_ps_1} yields outputs which are very close to the baseline Monte Carlo estimates, it is of interest to compare the computational cost of the two procedures. This is done in Table \ref{T_V_fixed}, which displays the mean run-time per discrete time step (first row) and the total run-time up to time $t_n$ (second row)
\begin{itemize}
\item for Algorithm \ref{alg_main_ps_1}, 
\item for the Monte Carlo method with $10^4$ samples and 
\item for a single run of the Euler-Maruyama scheme \eqref{sto_integrator}.
\end{itemize}
We see that the cost of running the moment-computing Algorithm \ref{alg_main_ps_1} is roughly of the same order as running the standard Euler-Maruyama scheme once, and three orders of magnitude less expensive than computing the Monte Carlo estimators.
}

\vspace{0cm}
\begin{table}[htb]
\cblack{
\begin{tabular}{|c||c|c|c|}
\hline
& Algorithm \ref{alg_main_ps_1} &Monte Carlo, 		&Euler-Maruyama,\\
&								&$10^4$ samples		&single run\\
\hline\hline
$\begin{matrix}
\text{Mean run-time}\\
\text{per time step}
\end{matrix}$ & $7.36\cdot 10^{-6}$\,s & $1.25\cdot 10^{-3}$\,s & $4.50\cdot 10^{-6}$\,s \\\hline
Total run-time & $\vphantom{A^{2^2}}1.16\cdot 10^1$\,s & $1.62\cdot 10^{3}$\,s & $5.84$\,s\\\hline
\end{tabular}
\caption{\cblack{Run-times in seconds (s) with fixed initial condition. The total number of discrete time steps is $n=1,296,000$.}}
\label{Table_Run_time_fix}
}
\end{table}

\begin{figure}[htb]
  \centering
  \hspace{0.cm}\includegraphics{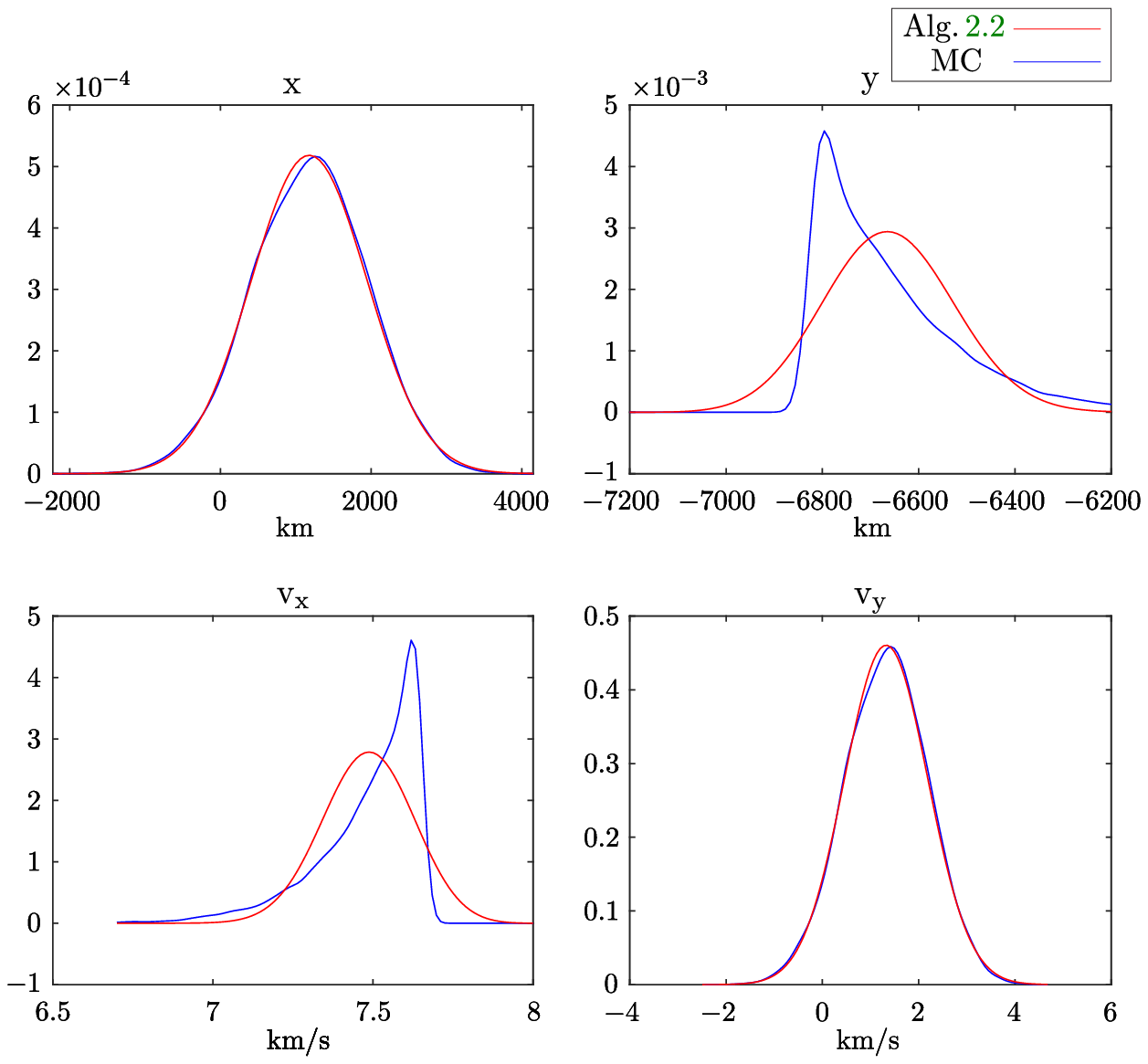}
  \caption{\cblack{Marginal pdf's at the final time $t_n$ when the initial condition is fixed. The red curves are the estimates obtained with Algorithm \ref{alg_main_ps_2} and the blue curves are KDEs computed from $10^4$ independent samples.}} 
  \label{Fig_1_fixed}
\end{figure}

\cblack{
Next, we turn attention to the performance of Algorithm \ref{alg_main_ps_2}, which yields estimates of the marginal densities of the state variables $x$, $y$, $v_x$ and $v_y$. Figure \ref{Fig_1_fixed} shows corresponding pdf's as generated by Algorithm \ref{alg_main_ps_2} (in red colour) and the kernel density estimators\footnote{We use the {\tt ksdensity} function available in Matlab, which determines the kernel bandwidth for the estimator automatically from the samples.} (KDEs) computed from the independent samples generated by running the Euler-Maruyama scheme \eqref{sto_integrator} $10^4$ times. We see that the KDEs are clearly non-Gaussian for $y$ and $v_x$, and the estimators computed via Algorithm \ref{alg_main_ps_2} fail to yield an accurate approximation in this case. Performance can be improved by increasing the order of the polynomial approximation, at the expense of a higher computational cost.
}

 
\cblack{Finally, Table \ref{T_V_fixed} shows the total variation distance (TVD) between the marginal densities estimated using Algorithm \ref{alg_main_ps_2} and Monte Carlo (i.e., a KDE with $10^4$ samples). Let us recall that the TVD between two probability distributions with pdf's $f$ and $g$ can be written as $\| f-g \|_{TV} = \frac{1}{2}\int_{\mathds{R}} |f(x)-g(x)|\diff x$. 
}

\begin{table}[htb]
\begin{center}
\begin{tabular}{|c|c|c|c|}\hline
$x$ & $y$ &$v_x$& $v_y$\\\hline\hline\cline{1-4}
$\vphantom{A^{2^2}}2.37\times 10^{-6}$&  $2.22\times 10^{-4}$ &  $1.84\times 10^{-1}$ &  $2.02\times 10^{-3}$\\\hline
\end{tabular}
\caption{\cblack{TVD between the estimates of the marginal densities computed with Algorithm \ref{alg_main_ps_2} and the KDEs computed from $10^4$ Monte Carlo samples.}}
\label{T_V_fixed}
\end{center}
\end{table}




\subsection{Random initial condition} \label{exper_2}

\cblack{
We illustrate the performance of the proposed approximation methods for the same dynamical model Section \ref{exper_1} except that the assume a random initial condition $\tb{X}_0$, modelled as a Gaussian random vector with mean $\tb{x}_0$ as in Eq. \eqref{in_cond_state} and covariance matrix
\begin{equation}\label{cov_mat_exp}
\tb{\Sigma}_0=\begin{pmatrix}
10^{-1} & 0 & 0 & 0\\
0 & 10^{-1} & 0 & 0\\
0 & 0 & 10^{-4} & 0\\
0 & 0 & 0 & 10^{-4}
\end{pmatrix}.
\end{equation}
The computation of moments is carried out using Algorithm \ref{alg_main_ps_1}, with a PCE of order $N_{\text{PCE}}=6$ and $N_s=420$ samples.
}

\begin{figure}[H]
\centering
	\includegraphics{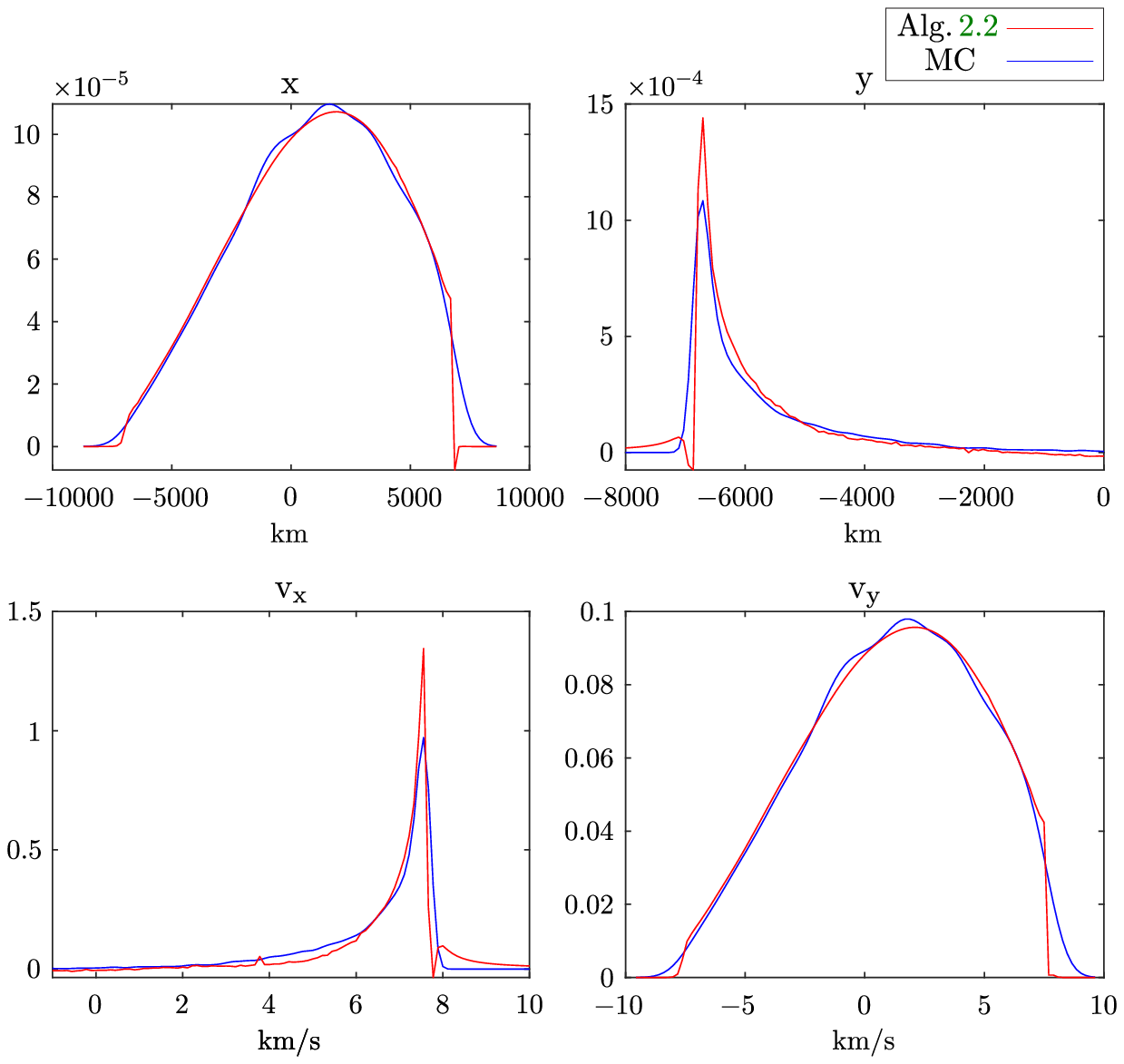}
    \caption{\cblack{Marginal pdf's at the final time $t_n$ when the initial condition is Gaussian. The red curves are the estimates obtained with Algorithm \ref{alg_main_ps_2} ($N_{PCE}=6$, $N_s=420$) and the blue curves are KDEs computed from $10^4$ independent samples.}} 
 \label{Fig_1_Un}
\end{figure}

\cblack{
The computer experiments are similar to Section \ref{exper_1}. In particular: 
\begin{itemize}
\item Table \ref{Table_mean_Un} shows a comparison of the expected values of the state variables ($x, y, v_x$ and $v_y$) as obtained though Algorithm \ref{alg_main_ps_1} and the baseline Monte Carlo method with $10^4$ independent trajectories. Both the expectations computed by the two methods and the absolute and relative differences are displayed. We observe small relative errors of order $10^{-3}$ for all state variables except $v_y$, which has a larger error of order $10^{-2}$.
\item Table \ref{Table_Cov_Un} shows a comparison of the second-order moment estimates, also at time $t_n$, using Algorithm \ref{alg_main_ps_1} and standard Monte Carlo with $10^4$ independent runs. The entry-wise relative differences between the two matrices, displayed in the fourth row of the table, shows nearly-uniform errors of order $10^{-2}$ for all variances and covariances.
\item Table \ref{Table_Run_time_Un} displays a comparison of the computational cost of Algorithm \ref{alg_main_ps_1} (with $N_{PCE}=6$ and $N_s=420$) and the baseline Monte Carlo procedure (with $10^4$ independendent runs). We see that the overall run-time of Algorithm \ref{alg_main_ps_1} is $\approx 67\%$ of the time needed for the Monte Carlo computation for a nearly identical performance. 
\item Figure \ref{Fig_1_Un} compares Algorithm \ref{alg_main_ps_2}, which estimates of the marginal pdf's of the state variables $x, y, v_x$ and $v_y$, and a KDE computed from the independent samples generated by running the Euler-Maruyama scheme \eqref{sto_integrator} $10^4$ times with random initialisations. All densities are clearly non-Gaussian and the two methods yield similar approximations. The performance of Algorithm \ref{alg_main_ps_2} can be improved by increasing the order of the polynomial approximation.
\item Finally, Table \ref{T_V_Un} shows the TVD between the marginal densities estimated using Algorithm \ref{alg_main_ps_2} and Monte Carlo (KDEs with $10^4$ samples). The largest error occurs for $v_x$, where the TVD is of order $10^{-2}$.
\end{itemize}
}

\vspace{1cm}\begin{table}[htbp]
\begin{center}
\begin{tabular}{|c|c|c|c|c|}\hline
&Algorithm \ref{alg_main_ps_1} 	&Monte Carlo,  	&Absolute	&Relative error\\
& 								&$10^4$ samples	&difference	&difference\\
\hline\hline
$x$ & $\vphantom{K^{B^2}}$ $\phantom{-}1.00\cdot 10^3$\,km & $\phantom{-}1.03\cdot 10^3$\,km & $2.71\cdot 10^{1}$\,km & $2.63\cdot 10^{-3}$\\\hline
$y$ & $\vphantom{K^{B^2}}-5.73\cdot 10^3$\,km & $-5.76\cdot 10^3$\,km & $3.10\cdot 10^{1}$\,km & $5.38\cdot 10^{-3}$\\\hline
$v_x$ & $\vphantom{K^{B^2}}6.43$\,km/s & $6.47$\,km/s & $3.50\cdot 10^{-2}$\,km/s & $5.40\cdot 10^{-3}$\\\hline
$v_y$ & $\vphantom{K^{B^2}_{B_2}}1.14$\,km/s & $1.17$\,km/s & $2.95\cdot 10^{-2}$\,km/s & $2.53\cdot 10^{-2}$\\\hline
\end{tabular}
\caption{\cblack{Estimate of $\mathds{E}[\tb{X}(t_n)]$ with the moment-computation Algorithm \ref{alg_main_ps_1}, compared with standard Monte Carlo estimates. The initial condition $\tb{X}_0$ is a Gaussian random vector.}}
\label{Table_mean_Un}
\end{center}
\end{table}

\begin{table}[H]
\begin{center}
\begin{tabular}{|c||c|}\hline
Algorithm \ref{alg_main_ps_1}&
$\vphantom{\begin{matrix}
\vspace{-0.3cm}a\\
a\\
a\\
a\\
a\\
a
\end{matrix}}\hspace{-0.cm}\begin{pmatrix}
\phantom{-}1.07\cdot 10^{7} &  \phantom{-} 1.56\cdot 10^{6} &-1.69\cdot 10^{3} & \!\!\!\!\phantom{-} 1.19\cdot 10^{4}\\\
\phantom{-}\!\!1.56\cdot 10^{6} & \phantom{-} 1.94\cdot 10^{6} &-2.17\cdot 10^{3} & \!\!\!\!\phantom{-} 1.75\cdot 10^{3} \\
 -1.69\cdot 10^{3} & -2.17\cdot 10^{3} &\phantom{-}\;\;\,2.44 & -1.90\\
\phantom{-}1.19\cdot 10^{4} & \phantom{-}1.75\cdot 10^{3} & \;\;\,-1.90 &\phantom{-}1.33\cdot 10^{1}
\end{pmatrix}$\\
\hline
Monte Carlo
&
$\vphantom{\begin{matrix}
\vspace{-0.3cm}a\\
a\\
a\\
a\\
a\\
a
\end{matrix}}\hspace{-0.cm}\begin{pmatrix}
 \phantom{-}1.03\cdot 10^{7} & \phantom{-}1.59\cdot 10^{6} & -1.73\cdot 10^{3} &  \!\!\!\!\phantom{-}1.15\cdot 10^{4}\\
 \phantom{-}1.59\cdot 10^{6} &  \phantom{-}1.89\cdot 10^{6} & -2.12\cdot 10^{3} &  \!\!\!\!\phantom{-}1.78\cdot 10^{3}\\
-1.73\cdot 10^{3} & -2.12\cdot 10^{3} & \;\;\,\phantom{-}2.38 & -1.94\\
\phantom{-}1.15\cdot 10^{4} & \phantom{-}1.78\cdot 10^{3}  & \;\;\,-1.94&\phantom{-}1.28\cdot 10^{1}
\end{pmatrix}$\hspace{-0.14cm}\\\hline
$\begin{matrix}
\text{Absolute}\\
\text{differences}
\end{matrix}
$
& 
$\vphantom{\begin{matrix}
\vspace{-0.3cm}a\\
a\\
a\\
a\\
a\\
a
\end{matrix}}\hspace{-0.1cm}
\begin{pmatrix}
 3.63\cdot 10^{5} &  2.75\cdot 10^{4} & \;\!\!\!\!\!3.54\cdot 10^{1} & 4.07\cdot 10^{2}\\
 2.75\cdot 10^{4} &  4.85\cdot 10^{4} & \;\!\!\!\!\!5.39\cdot 10^{1} & 3.02\cdot 10^{1}\\
 3.54\cdot 10^{1} & 5.39\cdot 10^{1} & \;6.01\cdot 10^{-2} & \;\;\,3.90\cdot 10^{-2}\\
4.07\cdot 10^{2} & 3.02\cdot 10^{1} & \;3.90\cdot 10^{-2} & \;\;\,4.56\cdot 10^{-1}
\end{pmatrix}$\hspace{-0.14cm}\\\hline
$\begin{matrix}
\text{Relative}\\
\text{differences}
\end{matrix}
$ & 
$\vphantom{\begin{matrix}
\vspace{-0.3cm}a\\
a\\
a\\
a\\
a\\
a
\end{matrix}}\hspace{-0.1cm}\begin{pmatrix}
3.52\cdot 10^{-2}& 1.73\cdot 10^{-2} & 2.05\cdot 10^{-2} & 3.54\cdot 10^{-2} \\
1.73\cdot 10^{-2}& 2.57\cdot 10^{-2} & 2.55\cdot 10^{-2} & 1.69\cdot 10^{-2} \\
2.05\cdot 10^{-2}& 2.55\cdot 10^{-2} & 2.53\cdot 10^{-2} & 2.01\cdot 10^{-2} \\
3.54\cdot 10^{-2}& 1.69\cdot 10^{-2} & 2.01\cdot 10^{-2} & 3.56\cdot 10^{-2} \\
\end{pmatrix}$\\\hline
\end{tabular}
\caption{\cblack{Estimates of the covariance matrix of $\tb{X}(t_n)$ computed via Algorithm \ref{alg_main_ps_1} and standard Monte Carlo, with $10^4$ independent samples. The initial condition $\tb{X}_0$ is random.}}
\label{Table_Cov_Un}
\end{center}
\end{table}

\vspace{0cm}
\begin{table}[H]
\begin{tabular}{|c||c|c|c|}
\hline
&Algorithm \ref{alg_main_ps_1} &Monte Carlo, 	&Euler-Maruyama, \\
&								&$10^4$ samples	&single run\\
\hline\hline
$\begin{matrix}
\text{Mean run-time}\\
\text{per step}
\end{matrix}$ & $8.38\times 10^{-4}$\,s & $1.25\times 10^{-3}$\,s & $4.55\times 10^{-6}$\,s \\\hline
\text{Total run-time} & $\vphantom{A^{2^2}}1.09\times 10^3$\,s & $1.62\times 10^{3}$\,s & $5.90$\,s\\\hline
\end{tabular}
\caption{\cblack{Run-times in seconds (s) with random initial condition ($N_{PCE}=6$, $N_s=420$). The total number of discrete time steps is $n=1,296,000$.}}
\label{Table_Run_time_Un}
\end{table}

\vspace{-0.8cm}

\vspace{-0.5cm}\begin{table}[H]
\begin{center}
\begin{tabular}{|c|c|c|c|}\hline
$x$ & $\textrm{y}$ &$v_x$& $v_y$\\
\hline\hline\cline{1-4}
$\vphantom{A^{2^2}}1.31\times 10^{-6}$&  $1.99\times 10^{-5}$ &  $1.72\times 10^{-2}$ &  $1.14\times 10^{-3}$\\\hline
\end{tabular}
\caption{\cblack{Total variation distance when initial condition is a Gaussian random vector. Algorithm \ref{alg_main_ps_2} has truncation order $N_{PCE}=6$ for the PCE and $N_s=420$ samples for the approximation. The Monte Carlo baseline KDEs are constructed from $10^4$ independent samples.}}
\label{T_V_Un}
\end{center}
\end{table}

%
%

%
\section{Conclusions} \label{sConclusions}

\cblack{
We introduced a methodology for the computation of the moments of the numerical solution of a multidimensional SDE, denoted $\hat{\tb{X}}_n$, using truncated Taylor polynomial approximations. The core of the method is the decomposition of the solution $\hat{\tb{X}}_n$ into a central part that can be computed deterministically from an ODE using an explicit numerical scheme and an effective noise process, whose moments determine the characterisation of $\hat{\tb{X}}_n$.
}

\cblack{
While we have derived the algorithm based on an Euler-Maruyama numerical scheme, the same ideas can be extended in a rather straightforward way to other explicit schemes, such as stochastic Runge-Kutta methods. When the initial condition is fixed, the proposed algorithm involves a single run of the Euler-Maruyama numerical scheme (plus some additional computations for the moments) and attains approximately the same performance as a Monte Carlos scheme with $10^4$ independent runs of the Euler-Maruyama scheme. When the initial condition is random, we resort to a PCE scheme and still attain the same performance as the standard Monte Carlo estimators of the mean and second order moments with just a fraction $\left( \approx \frac{2}{3} \right)$) of the run-time for a problem involving the propagation of uncertainty in a 2-dimensional Keplerian orbit. We have also shown how to use the approximate moments of the numerical solution to compute type A Gram-Charlier estimates of the 1-dimensional marginal pdf's of the dynamical variables. When the initial condition is random, the averaging due to the PCE scheme enables the approximation of densities which are clearly non-Gaussian.
}

\cblack{
The implementation of the algorithms as they have been presented demand the a priori calculation of the derivatives of the drift and diffusion coefficients. Although it has not been explored in this paper, such calculations can be implemented automatically in the numerical scheme resorting to the tools of Taylor differential algebra \cite{Valli13}.  
}


\bibliographystyle{siamplain}
\bibliography{references}
\end{document}

%% file: ex_shared.tex
\usepackage{lipsum}
\usepackage{amsfonts}
\usepackage{graphicx}
\usepackage{epstopdf}
\usepackage{algorithmic}
\ifpdf
  \DeclareGraphicsExtensions{.eps,.pdf,.png,.jpg}
\else
  \DeclareGraphicsExtensions{.eps}
\fi

\usepackage{savesym}
\usepackage{amssymb}
\usepackage{mathrsfs}
\usepackage{dsfont}
\usepackage{relsize}
\usepackage{array}
\usepackage{caption}
\usepackage{mathtools}
\usepackage{MnSymbol}
\savesymbol{dotplus}
\savesymbol{divideontimes}
\savesymbol{coloneq}
\savesymbol{corresponds}
\savesymbol{dagger}
\savesymbol{ddagger}
\savesymbol{mid}
\savesymbol{complement}
\savesymbol{barwedge}
\savesymbol{veebar}
\savesymbol{doublebarwedge}
\savesymbol{bigplus}
\savesymbol{bigtimes}
\savesymbol{bigcurlywedge}
\savesymbol{bigcurlyvee}
\savesymbol{bigominus}
\savesymbol{bigocirc}
\savesymbol{bigoslash}
\savesymbol{bigobackslash}
\usepackage{mathabx}
\usepackage{bbm}
\usepackage{fancyhdr,color}
\usepackage[bottom]{footmisc}
\usepackage{setspace}
\usepackage{hhline}
\usepackage{enumitem}
\usepackage{caption}

\newsiamremark{remark}{Remark}
\newsiamremark{hypothesis}{Hypothesis}
\crefname{hypothesis}{Hypothesis}{Hypotheses}
\newsiamthm{claim}{Claim}
\newsiamthm{assumption}{Assumption}

\newcommand{\diff}{\mathrm{d}}
\newcommand{\e}{\operatorname{e}}

\newcommand{\cblack}{\textcolor{black}}

\newcommand{\mF}{\mathcal{F}}
\newcommand{\mX}{\mathcal{X}}
\newcommand{\mbP}{\mathbb{P}}

\newcommand{\mbN}{\mathbb{N}}

\newcommand{\nn}{\nonumber}
\newcommand{\tb}{\boldsymbol}
\newcommand{\beq}{\begin{equation}}
\newcommand{\eeq}{\end{equation}}
\newcommand{\beqa}{\begin{eqnarray}}
\newcommand{\eeqa}{\end{eqnarray}}


\headers{Polynomial propagation of moments in stochastic differential...}{A. L\'opez Yela and J. M\'iguez}
\title{Polynomial propagation of moments in stochastic differential equations\thanks{This work has been partially supported by the European Space Agency (ESA contract No. 4000126151/19/D/SR, ``Uncertainty propagation meeting space debris needs(T711-501GR)'').}}

\author{Alberto L\'{o}pez Yela\thanks{Department of Signal Theory \& Communications, Universidad Carlos III de Madrid (Spain). E-mail: \email{alyela@math.uc3m.es}, \email{joaquin.miguez@uc3m.es}.}
\and Joaqu\'in M\'iguez\footnotemark[2]}

\usepackage{amsopn}
\DeclareMathOperator{\diag}{diag}
